\newcommand{\xbD}{\Delta}

\newcommand{\xbG}{\Gamma}

\newcommand{\xbO}{\Omega}
\newcommand{\xbP}{\Pi}

\newcommand{\xbS}{\Sigma}

\newcommand{\xba}{\alpha}
\newcommand{\xbb}{\beta}

\newcommand{\xbe}{\in}
\newcommand{\xbf}{\phi}
\newcommand{\xbg}{\gamma}

\newcommand{\xbm}{\mu}
\newcommand{\xbn}{\nu}
\newcommand{\xbo}{\omega}

\newcommand{\xbq}{\psi}

\newcommand{\xbs}{\sigma}

\newcommand{\xCK}{\times}

\newcommand{\xCN}{\neg}
\newcommand{\xCQ}{\emptyset}

\newcommand{\xCf}{\hspace{0.1em}}

\newcommand{\xcA}{\forall}

\newcommand{\xcC}{\not\subseteq}

\newcommand{\xcE}{\exists}

\newcommand{\xcL}{\not\vdash}
\newcommand{\xcM}{\not\models}
\newcommand{\xcN}{\hspace{0.2em}\not\sim\hspace{-0.9em}\mid\hspace{0.8em}}

\newcommand{\xcS}{\bigcap}
\newcommand{\xcT}{\bot}

\newcommand{\xcV}{\bigcup}

\newcommand{\xcX}{\Box}

\newcommand{\xcb}{\subset}
\newcommand{\xcc}{\subseteq}
\newcommand{\xcd}{\supseteq}
\newcommand{\xce}{\not\in}
\newcommand{\xcf}{\supset}
\newcommand{\xcg}{\geq}
\newcommand{\xch}{\Rightarrow}

\newcommand{\xcj}{\Leftrightarrow}
\newcommand{\xck}{\leq}
\newcommand{\xcl}{\vdash}
\newcommand{\xcm}{\models}
\newcommand{\xcn}{\hspace{0.2em}\sim\hspace{-0.9em}\mid\hspace{0.58em}}

\newcommand{\xco}{\vee}
\newcommand{\xcp}{\rightarrow}
\newcommand{\xcq}{\leftarrow}
\newcommand{\xcr}{\leftrightarrow}
\newcommand{\xcs}{\cap}
\newcommand{\xcu}{\wedge}
\newcommand{\xcv}{\cup}

\newcommand{\xcz}{\Box}

\newcommand{\xDC}{\hspace{2em}}
\newcommand{\xDH}{\item }

\newcommand{\xdA}{\mbox{\boldmath$A$}}
\newcommand{\xdB}{\mbox{\boldmath$B$}}
\newcommand{\xdC}{\mbox{\boldmath$C$}}
\newcommand{\xdD}{\mbox{\boldmath$D$}}

\newcommand{\xda}{{\cal A}}

\newcommand{\xdc}{{\cal C}}

\newcommand{\xdl}{{\cal L}}
\newcommand{\xdm}{{\cal M}}

\newcommand{\xdp}{{\cal P}}

\newcommand{\xdu}{{\cal U}}

\newcommand{\xdw}{{\cal W}}
\newcommand{\xdx}{{\cal X}}
\newcommand{\xdy}{{\cal Y}}
\newcommand{\xdz}{{\cal Z}}

\newcommand{\xEH}{ & }
\newcommand{\xEI}{\begin{itemize}}
\newcommand{\xEJ}{\end{itemize}}
\newcommand{\xEP}{ \\ }

\newcommand{\xEd}{\neq}
\newcommand{\xEh}{\begin{enumerate}}
\newcommand{\xEj}{\end{enumerate}}

\newcommand{\xeB}{\not\prec}

\newcommand{\xeb}{\prec}

\newcommand{\xee}{\succ}

\newcommand{\xej}{\lhd}

\newcommand{\xem}{\rhd}

\newcommand{\xer}{\sqsubset}

\newcommand{\xex}{\lceil}

\newcommand{\xFO}{\parallel}

\newcommand{\Xl}{\ldots}

\newcommand{\ol}{\overline}
\newcommand{\ul}{\underline}

\newcommand{\wt}{\overbrace}

\newcommand{\xssc}{\scriptsize}

\newcommand{\bl}{\begin{lemma} \rm}
\newcommand{\el}{\end{lemma}}
\newcommand{\br}{\begin{remark} \rm}
\newcommand{\er}{\end{remark}}
\newcommand{\be}{\begin{example} \rm}
\newcommand{\ee}{\end{example}}
\newcommand{\bco}{\begin{corollary} \rm}
\newcommand{\eco}{\end{corollary}}
\newcommand{\bc}{\begin{claim} \rm}
\newcommand{\ec}{\end{claim}}
\newcommand{\bfa}{\begin{fact} \rm}
\newcommand{\efa}{\end{fact}}
\newcommand{\bp}{\begin{proposition} \rm}
\newcommand{\ep}{\end{proposition}}
\newcommand{\bd}{\begin{definition} \rm}
\newcommand{\ed}{\end{definition}}
\newcommand{\bcs}{\begin{construction} \rm}
\newcommand{\ecs}{\end{construction}}
\newcommand{\bcd}{\begin{condition} \rm}
\newcommand{\ecd}{\end{condition}}
\newcommand{\bt}{\begin{theorem} \rm}
\newcommand{\et}{\end{theorem}}
\newcommand{\bn}{\begin{notation} \rm}
\newcommand{\en}{\end{notation}}
\newcommand{\bfi}{\begin{bild} \rm}
\newcommand{\efi}{\end{bild}}
\newcommand{\bsta}{\begin{statement} \rm}
\newcommand{\esta}{\end{statement}}
\newcommand{\bcom}{\begin{comment} \rm}
\newcommand{\ecom}{\end{comment}}
\newcommand{\bdia}{\begin{diagram} \rm}
\newcommand{\edia}{\end{diagram}}

\newcommand{\bfc}{\begin{figure}[htb] \begin{center}}
\newcommand{\efc}{\end{center} \end{figure}}

\documentclass{article}
\usepackage{amssymb}

\usepackage{epic,eepic}

\title{
A THEORY OF HIERARCHICAL CONSEQUENCE AND CONDITIONALS
}
\author{Dov M Gabbay
\thanks{
Dov.Gabbay@kcl.ac.uk, www.dcs.kcl.ac.uk/staff/dg
} \\
King's College, London
\thanks{
Department of Computer Science, King's College London, Strand,
London WC2R 2LS, UK
} \\ \\
Karl Schlechta
\thanks{
ks@cmi.univ-mrs.fr, karl.schlechta@web.de, http://www.cmi.univ-mrs.fr/ $\sim$ ks
} \\
Laboratoire d'Informatique Fondamentale de Marseille
\thanks{
UMR 6166, CNRS and Universit\'{e} de Provence,
Address: CMI, 39, rue Joliot-Curie, F-13453 Marseille Cedex 13, France
}
}

\begin{document}

\newtheorem{lemma}{Lemma}[section]
\newtheorem{theorem}[lemma]{Theorem}
\newtheorem{proposition}[lemma]{Proposition}
\newtheorem{corollary}[lemma]{Corollary}
\newtheorem{claim}[lemma]{Claim}
\newtheorem{fact}[lemma]{Fact}
\newtheorem{remark}[lemma]{Remark}
\newtheorem{definition}{Definition}[section]
\newtheorem{construction}{Construction}[section]
\newtheorem{condition}{Condition}[section]
\newtheorem{example}{Example}[section]
\newtheorem{notation}{Notation}[section]
\newtheorem{bild}{Figure}[section]
\newtheorem{comment}{Comment}[section]
\newtheorem{statement}{Statement}[section]
\newtheorem{diagram}{Diagram}[section]

\maketitle

\renewcommand{\labelenumi}
  {(\arabic{enumi})}
\renewcommand{\labelenumii}
  {(\arabic{enumi}.\arabic{enumii})}
\renewcommand{\labelenumiii}
  {(\arabic{enumi}.\arabic{enumii}.\arabic{enumiii})}
\renewcommand{\labelenumiv}
  {(\arabic{enumi}.\arabic{enumii}.\arabic{enumiii}.\arabic{enumiv})}

\setcounter{secnumdepth}{3}
\setcounter{tocdepth}{3}

\begin{abstract}

We introduce $\xda$-ranked preferential structures and combine them with an
accessibility relation. $\xda$-ranked preferential structures are intermediate
between simple preferential structures and ranked structures. The additional
accessibilty relation allows us to consider only parts of the overall
$\xda$-ranked structure. This framework allows us to formalize contrary to duty
obligations, and other pictures where we have a hierarchy of situations,
and maybe not all are accessible to all possible worlds. Representation
results are proved.

\end{abstract}

\tableofcontents

% ******* BEGIN LATEX SOURCE FILE hier-1.tex *******
%
% Uebers. aus Karltex File: hier-1.ms
%
%
%  A THEORY OF HIERARCHICAL CONDITIONALS (Paper 319)
%  A THEORY OF HIERARCHICAL CONDITIONALS (Paper 319)
% %
% ==================================================
%  INTRODUCTION
%  INTRODUCTION
% %
% =============
\section{
Introduction
}

% {\LARGE karl-search= Start Hier-In }

\label{Section Hier-In}
\index{Section Hier-In}
\subsection{
Description of the problem
}

This paper, like all papers about nonmonotonic logics, is about
formalization
of (an aspect of) common sense reasoning.

We often see a hierarchy of situations, e.g.:

 \xEh

 \xDH it is better to prevent an accident than to help the victims,

 \xDH it is better to prove a difficult theorem than to prove an easy
lemma,

 \xDH it is best not to steal, but if we have stolen, we should return
the stolen object to its legal owner, etc.

 \xEj

On the other hand, it is sometimes impossible to achieve the best
objective.

We might have seen the accident happen from far away, so we were unable to
interfere in time to prevent it, but we can still run to the scene and
help
the victims.

We might have seen friends last night and had a drink too many, so today's
headaches will not allow us to do serious work, but we can still prove a
little lemma.

We might have needed a hammer to smash the windows of a car involved in an
accident, so we stole it from a building site, but will return it
afterwards.

We see in all cases:

- a hierarchy of situations

- not all situations are possible or accessible for an agent.

In addition, we often have implicitly a ``normality'' relation:

Normally, we should help the victims, but there might be situations where
not: This would expose ourselves to a very big danger, or this would
involve
neglecting another, even more important task (we are supervisor in a
nuclear
power plant  \Xl.), etc.

Thus, in all ``normal'' situations where an accident seems imminent, we
should
try to prevent it. If this is impossible, in all ``normal'' situations, we
should
help the victims, etc.

We combine these three ideas

 \xEh
 \xDH normality,
 \xDH hierarchy,
 \xDH accessibility
 \xEj

in the present paper.

Note that it might be well possible to give each situation a numerical
value and decide by this value what is right to do - but humans do not
seem
to think this way, and we want to formalize human common sense reasoning.

Before we begin the formal part, we elaborate above situations with more
examples.

 \xEI

 \xDH

We might have the overall intention to advance computer science.

So we apply for the job of head of department of computer science at
Stanford,
and promise every tenured scientist his own laptop.

Unfortunately, we do not get the job, but become head of computer
science department at the local community college. The college does not
have research as priority, but we can still do our best to achieve our
overall intention, by, say buying good books for the library, or buy
computers
for those still active in research, etc.

So, it is reasonable to say that, even if we failed in the best possible
situation - it was not accessible to us - we still succeeded in another
situation, so we achieved the overall goal.

 \xDH

The converse is also possible, where better solutions become possible,
as is illustrated by the following example.

The daughter and her husband say to have the overall intention to start a
family
life with a house of their own, and children.

Suppose the mother now asks her daughter: You have been married now for
two
years, how come you are not pregnant?

Daughter - we cannot afford a baby now, we had to take a huge mortgage to
buy
our house and we both have to work.

Mother - $I$ shall pay off your mortgage. Get on with it!

In this case, what was formerly inaccessible, is now accessible, and
if the daughter was serious about her intentions - the mother can begin
to look for baby carriages.

Note that we do not distinguish here how the situations change, whether
by our own doing, or by someone else's doing, or by some events not
controlled
by anyone.

 \xDH

Consider the following hierarchy of obligations making fences as
unobtrusive as possible, involving contrary to duty obligations.

 \xEh

 \xDH You should have no fence (main duty).

 \xDH If this is impossible (e.g. you have a dog which might invade
neighbours'
property), it should be less than 3 feet high (contrary to duty, but
second
best choice).

 \xDH If this is impossible too (e.g. your dog might jump over it), it
should be
white (even more contrary to duty, but still better than nothing).

 \xDH If all is impossible, you should get the neighbours' consent (etc.).

 \xEj

 \xEJ
\subsection{
Outline of the solution
}

The last example can be modelled as follows $( \xbm (x)$ is the minimal
models of $x):$

Layer 1: $ \xbm (True):$ all best models have no fence.

Layer 2: $ \xbm (fence):$ all best models with a
fence are less than 3 ft. high.

Layer 3: $ \xbm (fence$ and more than 3 ft. high): all best models with a
tall
fence have a white fence.

Layer 4: $ \xbm (fence$ and non-white and $ \xcg 3$ ft): in all best
models with a
non-white fence taller than 3 feet, you have permission

Layer 5: all the rest

This will be modelled by a corresponding $ \xda -$structure.

In summary:

 \xEh

 \xDH We have a hierarchy of situations, where one group (e.g. preventing
accidents) is strictly better than another group (e.g. helping victims).

 \xDH Within each group, preferences are not so clear (first help person
A,
or person $B,$ first call ambulance, etc.?).

 \xDH We have a subset of situations which are attainable, this can be
modelled
by an accessibility relation which tells us which situations are possible
or
can be reached.

 \xEj

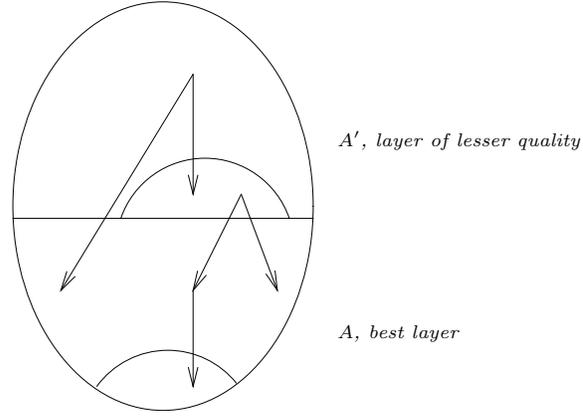
\vspace{20mm}

\begin{diagram}

\label{Diagram A-Ranked}
\index{Diagram A-Ranked}

\centering
\setlength{\unitlength}{0.00083333in}
{\renewcommand{\dashlinestretch}{30}
\begin{picture}(2390,3581)(0,0)
\put(1212.000,2028.000){\arc{1110.000}{3.4719}{5.9529}}
\put(979.651,835.201){\arc{1095.700}{3.7717}{5.6071}}
\put(949,2283){\ellipse{1874}{2550}}
\path(12,2208)(1887,2208)
\path(1137,3108)(1137,2358)
\path(1107.000,2478.000)(1137.000,2358.000)(1167.000,2478.000)
\path(1137,3108)(312,1758)
\path(348.976,1876.037)(312.000,1758.000)(400.172,1844.750)
\path(1437,2358)(1662,1758)
\path(1591.775,1859.826)(1662.000,1758.000)(1647.955,1880.893)
\path(1437,2358)(1137,1758)
\path(1163.833,1878.748)(1137.000,1758.000)(1217.498,1851.915)
\path(1137,1758)(1137,1158)
\path(1107.000,1278.000)(1137.000,1158.000)(1167.000,1278.000)
\put(2037,2658){{\xssc $A'$, layer of lesser quality}}
\put(2037,1458){{\xssc $A$, best layer}}
\put(-700,800){{\xssc
Each layer behaves inside like any preferential structure.}}
\put(-700,600){{\xssc
Amongst each other, layers behave like ranked structures.}}

\put(100,200) {{\rm\bf $\xda-$ ranked structure}}

\end{picture}
}
\end{diagram}

\vspace{4mm}

We combine all three ideas, consider what we call $ \xda -$ranked
structures,
structures which are organized in levels $A_{1},$ $A_{2},$ $A_{3},$ etc.,
where all
elements of $A_{1}$ are better than any element of $A_{2}$ - this is
basically
rankedness -, and where inside each $A_{i}$ we have an arbitrary relation
of
preference. Thus, an $ \xda -$ranked structure is between a simple
preferential
structure and a fully ranked structure.

See Diagram \ref{Diagram A-Ranked}.

Remark: It is not at all necessary that the rankedness relation between
the different layers and the relation inside the layers express the
same concept. For instance, rankedness may express deontic preference,
whereas the inside relation expresses normality or some usualness.

In addition, we have an accessibility relation $R,$ which tells us which
situations are reachable.

It is perhaps easiest to motivate the precise choice of modelling
by layered (or contrary to duty) obligations.

For any point $t,$ let $R(t):=\{s:tRs\},$ the set of $R-$reachable points
from $t.$ Given a preferential structure $ \xdx:=<X, \xeb >,$ we can
relativize $ \xdx $
by considering only those points in $X,$ which are reachable from $t.$

Let $X' \xcc X,$ and $ \xbm (X' )$ the minimal points of $X,$ we will now
consider
$ \xbm (X' ) \xcs R(t)$ - attention, not: $ \xbm (X' \xcs R(t))!$ This
choice is motivated
by the following: norms are universal, and do not depend on one's
situation $t.$

If $ \xdx $ describes a simple obligation, then we are obliged to $Y$ iff
$ \xbm (X' ) \xcs R(t) \xEd \xCQ,$ and $ \xbm (X' ) \xcs R(t) \xcc Y.$
The first clause excludes
obligations to the unattainable. We can write this as follows, supposing
that $X' $ is the set of models of $ \xbf ',$ and $Y$ is the set of
models of $ \xbq:$

$m \xcm \xbf ' > \xbq.$

Thus, we put the usual consequence relation $ \xcn $ into the object
language
as $>,$ and relativize to the attainable (from $m).$

If an $ \xda -$ranked structure has two or more layers, then we are, if
possible,
obliged to fulfill the lower obligation, e.g. prevent an accident,
but if this is impossible, we are obliged to fulfill the upper
obligation, e.g. help the victims, etc.

See Diagram \ref{Diagram Pischinger}.

\vspace{10mm}

\begin{diagram}

\label{Diagram Pischinger}
\index{Diagram Pischinger}

\centering
\setlength{\unitlength}{1mm}
{\renewcommand{\dashlinestretch}{30}
\begin{picture}(150,150)(0,0)

\path(30,10)(30,110)(110,110)(110,10)(30,10)
\path(30,30)(110,30)
\path(30,50)(110,50)
\path(30,70)(110,70)
\path(30,90)(110,90)

\put(70,10){\arc{20}{-3.14}{0}}
\put(70,30){\arc{20}{-3.14}{0}}
\put(70,50){\arc{20}{-3.14}{0}}
\put(70,70){\arc{20}{-3.14}{0}}
\put(70,90){\arc{20}{-3.14}{0}}

\put(60,80){\circle{30}}

\path(5,80)(25,80)
\path(22.3,81)(25,80)(22.3,79)
\put(5,80){\circle*{1}}
\put(25,80){\circle*{1}}

\put(15,75){\xssc{$R$}}
\put(5,75){\xssc{$t$}}
\put(25,75){\xssc{$s$}}

\put(10,130){\xssc{The overall structure is visible from $t$}}
\put(10,123){\xssc{Only the inside of the circle is visible from $s$}}
\put(10,116){\xssc{Half-circles are the sets of minimal elements of layers}}

\put(20,3) {{\rm\bf $\xda-$ ranked structure and accessibility}}

\end{picture}
}

\end{diagram}

\vspace{4mm}

Let now, for simplicity, $ \xdB $ be a subset of the union of all layers
A, and
let $ \xdB $ be the set of models of $ \xbb.$ This can be done, as the
individual
subset can be found by considering $A \xcs \xdB,$ and call the whole
structure
$< \xda, \xdB >.$

Then we say that $m$ satisfies $< \xda, \xdB >$ iff in the lowest layer A
where
$ \xbm (A) \xcs R(m) \xEd \xCQ $ $ \xbm (A) \xcs R(m) \xcc \xdB.$

When we want a terminology closer to usual conditionals, we may
write e.g. $(A_{1}>B_{1};A_{2}>B_{2}; \Xl.)$ expressing that the best is
$A_{1},$ and
then $B_{1}$ should hold, the second best is $A_{2},$ then $B_{2}$ should
hold, etc.
(The $B_{i}$ are just $A_{i} \xcs \xdB.)$
See Diagram \ref{Diagram C-Validity}.
\subsection{
Historical remarks
}

 \xEh

 \xDH In an abstract consideration of desirable properties a logic might
have,
 \cite{Gab85}
examined rules a nonmonotonic consequence relation $ \xcn $ should
satisfy:

 \xEh

 \xDH (REF) $ \xbD, \xba \xcn \xba,$

 \xDH (CUM) $ \xbD \xcn \xba $ $ \xch $ $( \xbD \xcn \xbb $ $ \xcj $ $
\xbD, \xba \xcn \xbb ).$

 \xEj

Preferential structures themselves were introduced as abstractions
of Circumscription independently in  \cite{Sho87b} and  \cite{BS85}.
A precise definition of these structures is given below
in Definition \ref{Definition Pref-Str}.

Both, the semantic and the syntactic, approaches were connected
in  \cite{KLM90}, where a representation
theorem was proved, showing that the (stronger than $Gabbay' $s) system
$P$
corresponds to ``smooth'' preferential structures. System $P$ consists of

 \xEh

 \xDH (AND) $ \xbf \xcn \xbq,$ $ \xbf \xcn \xbq ' $ $ \xch $ $ \xbf \xcn
\xbq \xcu \xbq ',$
 \xDH (OR) $ \xbf \xcn \xbq,$ $ \xbf ' \xcn \xbq $ $ \xch $ $ \xbf \xco
\xbf ' \xcn \xbq,$
 \xDH (LLE) $ \xcl \xbf \xcr \xbf ' $ $ \xch $ $( \xbf \xcn \xbq $ $ \xcj
$ $ \xbf ' \xcn \xbq ),$
 \xDH (RW) $ \xbf \xcn \xbq,$ $ \xcl \xbq \xcp \xbq ' $ $ \xch $ $ \xbf
\xcn \xbq ',$
 \xDH (SC) $ \xcl \xbf \xcp \xbf ' $ $ \xch $ $ \xbf \xcn \xbf ',$
 \xDH (CUM) $ \xbf \xcn \xbq $ $ \xch $ $( \xbf \xcn \xbq ' $ $ \xcj $ $
\xbf \xcu \xbq \xcn \xbq ' ).$

 \xEj

where $ \xcl $ is classical provability.

Details can be found in Definition \ref{Definition Log-Cond}.

 \xDH Ranked preferential structures were introduced in  \cite{LM92},
see Definition \ref{Definition Rank-Rel}. On the logical side, they
correspond to above
system $P,$ plus the additional axiom:

(RatM) $ \xbf \xcn \xbq,$ $ \xbf \xcN \xCN \xbq ' $ $ \xch $ $ \xbf \xcu
\xbq ' \xcn \xbq.$

 \xDH Accessibility relations in possible worlds semantics go back (at
least) to
Kripke's semantics for modal logics.

 \xEj
\subsection{
Formal modelling and summary of results
}

% {\LARGE karl-search= Start Hier-In-Summary }

\label{Section Hier-In-Summary}
\index{Section Hier-In-Summary}

We started with an investigation of ``best fulfillment'' of abstract
requirements,
and contrary to duty obligations.

It soon became evident that semi-ranked preferential structures give a
natural
semantics to contrary to duty obligations, just as simple preferential
structures give a natural semantics to simple obligations - the latter
goes
back to Hansson  \cite{Han69}.

A semi-ranked - or $ \xda -$ranked preferential
structure, as we will call them later, as they are based on a system of
sets $ \xda $ - has a finite number of layers, which amongst them are
totally ordered
by a ranking, but the internal ordering is just any (binary) relation.
It thus has stronger properties than a simple preferential structure, but
not as strong ones as a (totally) ranked structure.

The idea is to put the (cases of the) strongest obligation at
the bottom, and the weaker ones more towards the top. Then, fulfillment of
a strong obligation makes the whole obligation automatically satisfied,
and
the weaker ones are forgotten.

Beyond giving a natural semantics to contrary to duty obligations,
semi-ranked
structures seem very useful for other questions of knowledge
representation.
For instance, any blackbird might seem a more normal bird than any
penguin,
but we might not be so sure within each set of birds.

Thus, this generalization of preferential semantics seems very natural and
welcome.

The second point of this paper is to make some, but not necessarily all,
situations accessible to each point of departure. Thus, if we imagine
agent $ \xCf a$
to be at point $p,$ some fulfillments of the obligation, which are
reachable
to agent $a' $ from point $p' $ might just be impossible to reach for him.
Thus, we introduce a second relation, of accessibility in the intuitive
sense, denoting situations which can be reached. If this relation is
transitive, then we have restrictions on the set of reachable situations:
if $p$ is accessible from $p',$ and $p$ can access situation $s,$ then so
can
$p',$ but not necessarily the other way round.

On the formal side, we characterize:

(1) $ \xda -$ranked structures,

(2) satisfaction of an $ \xda -$ranked conditional once an accessibility
relation
between the points $p,$ $p',$ etc. is given.

For the convience of the reader, we now state the main formal results of
this paper - together with the more unusual definitions.

On (1):

Let $ \xdA $ be a fixed set, and $ \xda $ a finite, totally ordered (by
$<)$ disjoint cover
by non-empty subsets of $ \xdA.$

For $x \xbe \xdA,$ let $rg(x)$ the unique $A \xbe \xda $ such that $x
\xbe A,$ so $rg(x)<rg(y)$ is
defined in the natural way.

A preferential structure $< \xdx, \xeb >$ $( \xdx $ a set of pairs
$<x,i>)$ is called $ \xda -$ranked
iff for all $x,x' $ $rg(x)<rg(x' )$ implies $<x,i> \xeb <x',i' >$ for all
$<x,i>,<x',i' > \xbe \xdx.$
See Definition \ref{Definition Pref-Str} for the definition of preferential
structures, and Diagram \ref{Diagram A-Ranked} for an illustration.

We then have:

Let $ \xcn $ be a logic for $ \xdl.$ Set $T^{ \xdm }:=Th( \xbm_{ \xdm
}(M(T))),$ and $ \ol{ \ol{T} }:=\{ \xbf:T \xcn \xbf \}.$
where $ \xdm $ is a preferential structure.

(1) Then there is a (transitive) definability preserving
classical preferential model $ \xdm $ s.t. $ \ol{ \ol{T} }=T^{ \xdm }$ iff

(LLE), (CCL), (SC), (PR) hold for all $T,T' \xcc \xdl.$

(2) The structure can be chosen smooth, iff, in addition

(CUM) holds.

(3) The structure can be chosen $ \xda -$ranked, iff, in addition

$( \xda -$min) $T \xcL \xCN \xba_{i}$ and $T \xcL \xCN \xba_{j},$ $i<j$
implies $ \ol{ \ol{T} } \xcl \xCN \xba_{j}$

holds.

See Definition \ref{Definition Pref-Log} for the logic defined by a
preferential
structure, Definition \ref{Definition Log-Cond} for the logical conditions,
Definition \ref{Definition Smooth} for smoothness.

On (2)

Given a transitive accessibility relation $R,$ $R(m):=\{x:mRx\}.$

Given $ \xda $ as above, let $ \xdB \xcc \xdA $ be the set of ``good''
points in $ \xdA,$ and
set $ \xdc:=< \xda, \xdB >.$

We define:

(1) $ \xbm ( \xda ):= \xcV \{ \xbm (A_{i}):i \xbe I\}$

(warning: this is NOT $ \xbm ( \xdA ))$

(2) $ \xda_{m}:=R(m) \xcs \xdA,$

(3) $ \xbm ( \xda_{m}):= \xcV \{ \xbm (A_{i}) \xcs R(m):i \xbe I\}$

(3a) $ \xbn ( \xda_{m})$ $:=$ $ \xbm ( \xbm ( \xda_{m}))$

(thus $ \xbn ( \xda_{m})$ $=$ $\{a \xbe \xdA:$ $ \xcE A \xbe \xda (a \xbe
\xbm (A),$ $a \xbe R(m),$ and

$ \xDC \xDC \xDC \xCN \xcE a' ( \xcE A' \xbe \xda (a' \xbe \xbm (A' ),$
$a' \xbe R(m),$ $a' \xeb a\}.$

(4) $m \xcm \xdc $ $: \xcr $ $ \xbn ( \xda_{m})) \xcc \xdB.$

See Diagram \ref{Diagram C-Validity}

Then the following hold:

Let $m,m' \xbe M,$ $A,A' \xbe \xda,$ $ \xdA $ be the set of models of $
\xba.$

(1) $m \xcm \xcX \xCN \xba,$ mRm' $ \xch $ $m' \xcm \xcX \xCN \xba $

(2) $ \xCf mRm',$ $ \xbn ( \xda_{m}) \xcs A \xEd \xCQ,$ $ \xbn (
\xda_{m' }) \xcs A' \xEd \xCQ,$ $ \xch $ $A \xck A' $ (in the ranking)

(3) $ \xCf mRm',$ $ \xbn ( \xda_{m}) \xcs A \xEd \xCQ,$ $ \xbn (
\xda_{m' }) \xcs A' \xEd \xCQ,$ $m \xcm \xdc,$ $m' \xcM \xdc, \xch $
$A<A' $

Conversely, these conditions suffice to construct an accessibility
relation between $M$ and $ \xdA $ satisfying them, so they are sound and
complete.

% karl-search= End Hier-In-DovIn-Motiv
\vspace{7mm}

% *************************************

\vspace{7mm}

% karl-search= End Hier-In-DovIn
\vspace{7mm}

% *************************************

\vspace{7mm}

% karl-search= End Hier-In-Summary
\vspace{7mm}

% *************************************

\vspace{7mm}

\subsection{
Overview
}

We next point out some connections with other domains of artificial
intelligence and computer science.

We then put our work in perspective with a summary of logical and
semantical conditions for nonmonotonic and related logics, and present
basic defintions for preferential structures.

Next, we will give special definitions for our framework.

We then start the main formal part, and prove representation results
for $ \xda -$ranked structures, first for the general case, then for
the smooth case. The general case needs more work, as we have to
do a (minor) modification of the not $ \xda -$ranked case. The smooth case
is
easy, we simply have to append a small construction. Both proofs are
given in full detail, in order to make the text self-contained.

Finally, we characterize changes due to restricted accessibility.

% karl-search= End Hier-In
\vspace{7mm}

% *************************************

\vspace{7mm}

\section{
Connections with other concepts
}
\subsection{
Hierarchical conditionals and programs
}

% {\LARGE karl-search= Start Hier-Cond }

\label{Section Hier-Cond}
\index{Section Hier-Cond}

Our situation is now very similar to a sequence of computer program
instructions:

if $A_{1}$ then do $B_{1};$

else if $A_{2}$ then do $B_{2};$

else if $A_{3}$ then do $B_{3};$

where we can see the $B_{i}$ as subroutines.

We can deepen this analogy in two directions:

(1) connect it to Update

(2) put an imperative touch to it.

In both cases, we differentiate between different degrees of fulfillment
of $ \xdc:$ the lower the level is which is fulfilled, the better.

(1) We can consider all threads of reachability which lead to a model $m$
where
$m \xcm \xdc.$ Then we take as best threads those which lead to the best
fulfillment of
$ \xdc.$ So degree of fulfillment gives the order by which we should do
the update.
(This is then not update in the sense that we choose the most normal
developments, but rather we actively decide for the most desirable ones.)
We will not pursue this line any further here, but leave it for future
research.

(2): We introduce an imperative operator, say!.! means that one should
fulfill $ \xdc $ as best as possible by suitable choices.
We will elaborate this now.

First, we can easily compare the degree of satisfaction of $ \xdc $ of two
models:

\bd

$\hspace{0.01em}$

% (+++*** Orig. No.:  Definition 3.1 )

\label{Definition 3.1}

Let $m,m' \xcm \xdc,$ and define $m<m' $ $: \xch $ $ \xbm ( \xbm (
\xda_{m}) \xcv \xbm ( \xda_{m' })) \xcs \xbm ( \xda_{m' })= \xCQ.$ $(
\xbm $ is,
as usual, relative to some fixed $ \xck_{t}.)$

\ed

For two sets of models, $X,$ $X',$ the situation does not seem so easy.
So suppose
that $X,X' \xcm \xdc.$ First, we
have to decide how to compare this, we do by the maximum: $X<X' $ iff the
worst
satisfaction of all $x \xbe X$ is better than the worst satisfaction in
$X'.$
More precisely, we look at all $ \xbg ( \xdc )$ for all $x \xbe X,$ take
the maximum (which
exists, as $ \xda $ is finite), and then compare the maxima for $X$ and
for $X'.$

Suppose now that there are points where we can make decisions $('' free$
$will'' ),$
let $m$ be such a point. We introduce a new relation $D,$ and let mDm' iff
we can
decide to go from $m$ to $m'.$ The relation $D$ expresses this
possibility - it is
our definition of ``free will''.

\bd

$\hspace{0.01em}$

% (+++*** Orig. No.:  Definition 3.2 )

\label{Definition 3.2}

Consider now some formula $ \xbf,$ and define

$m \xcm! \xbf $ $: \xch $ $D(m) \xcs M( \xbf )<D(m) \xcs M( \xCN \xbf )$

(as defined in Definition \ref{Definition 3.1}).

% karl-search= End Hier-Cond
\vspace{7mm}

% *************************************

\vspace{7mm}

\subsection{
Connection with Theory Revision
}

% {\LARGE karl-search= Start Hier-TR }

\label{Section Hier-TR}
\index{Section Hier-TR}

\ed

In particular, the situation of contrary to duty obligations
(see Section \ref{Section Hier-Oblig}) shows an intuitive
similarity to revision. You have the duty not to have a fence. If this is
impossible (read: inconsistent), then it should be white. So the duty is
revised.

But there is also a formal analogy: As is well known, AGM revision (with
fixed left hand side $K)$ corresponds to a ranked order of models, where
models of $K$ have lowest rank (or: distance 0 from $K-$models). The
structures we consider $( \xda -$rankings) are partially ranked, i.e.
there is only a partial ranked preference, inside the layers, nothing is
said about the ordering. This partial ranking is natural, as we have only
a limited number of cases to consider.

But we use the revision order (based on $K,$ so it really is a $ \xck_{K}$
relation)
differently: We do not revise $K,$ but use only the order to choose the
first
layer which has non-empty intersection with the set of possible cases.
Still, the spirit (and formal apparatus) of revision is there, just used
somewhat differently. The $K-$relation expresses here deontic quality, and
if
the best situation is impossible, we choose the second best, etc.

Theory revision with variable $K$ is expressed by a distance between
models
(see  \cite{LMS01}),
where $K* \xbf $ is defined by the set of $ \xbf $ models which have
minimal distance from
the set of $K$ models.

We can now generalize our idea of layered structure to a partial distance
as
follows: For instance, $d(K,A)$ is defined, $d(K,B)$ too, and we know that
all A models with minimal distance to $K$ have smaller distance than the
$B$ models
with minimal distance to $K.$ But we do NOT know a precise distance for
other A
models, we can sometimes compare, but not always. We may also know that
all
A models are closer to $K$ than any $B$ model is, but for a and $a',$
both A models,
we might not know if one or the other is closer to $K,$ or is they have
the
same distance.

% karl-search= End Hier-TR
\vspace{7mm}

% *************************************

\vspace{7mm}

\subsection{
Definitions for our framework
}
%  6  REPRESENTATION for &a-ranked structures
%  6  REPRESENTATION for &a-ranked structures
% %
% ===========================================
\section{
Representation results for $\xda$-ranked structures
}

% {\LARGE karl-search= Start Hier-ARepr }

\label{Section Hier-ARepr}
\index{Section Hier-ARepr}
%  6.1  Basic definitions for preferential structures
%  6.1  Basic definitions for preferential structures
% %
% ===================================================
\subsection{
Basic definitions for preferential structures
}

% {\LARGE karl-search= Start Hier-ARepr-Def }

\label{Section Hier-ARepr-Def}
\index{Section Hier-ARepr-Def}
\index{Definition Alg-Base}

\bd

$\hspace{0.01em}$

% (+++*** Orig. No.:  Definition Alg-Base )

\label{Definition Alg-Base}

We use $ \xdp $ to denote the power set operator,
$ \xbP \{X_{i}:i \xbe I\}$ $:=$ $\{g:$ $g:I \xcp \xcV \{X_{i}:i \xbe I\},$
$ \xcA i \xbe I.g(i) \xbe X_{i}\}$ is the general cartesian
product, $card(X)$ shall denote the cardinality of $X,$ and $V$ the
set-theoretic
universe we work in - the class of all sets. Given a set of pairs $ \xdx
,$ and a
set $X,$ we denote by $ \xdx \xex X:=\{<x,i> \xbe \xdx:x \xbe X\}.$ When
the context is clear, we
will sometime simply write $X$ for $ \xdx \xex X.$

$A \xcc B$ will denote that $ \xCf A$ is a subset of $B$ or equal to $B,$
and $A \xcb B$ that $ \xCf A$ is
a proper subset of $B,$ likewise for $A \xcd B$ and $A \xcf B.$

Given some fixed set $U$ we work in, and $X \xcc U,$ then $ \xdC (X):=U-X$
.

If $ \xdy \xcc \xdp (X)$ for some
$X,$ we say that $ \xdy $ satisfies

$( \xcs )$ iff it is closed under finite intersections,

$( \xcS )$ iff it is closed under arbitrary intersections,

$( \xcv )$ iff it is closed under finite unions,

$( \xcV )$ iff it is closed under arbitrary unions,

$( \xdC )$ iff it is closed under complementation.

We will sometimes write $A=B \xFO C$ for: $A=B,$ or $A=C,$ or $A=B \xcv
C.$

We make ample and tacit use of the Axiom of Choice.
\index{Definition Rel-Base}

\ed

\bd

$\hspace{0.01em}$

% (+++*** Orig. No.:  Definition Rel-Base )

\label{Definition Rel-Base}

$ \xeb^{*}$ will denote the transitive closure of the relation $ \xeb.$
If a relation $<,$
$ \xeb,$ or similar is given, $a \xcT b$ will express that a and $b$ are
$<-$ (or $ \xeb -)$
incomparable - context will tell. Given any relation $<,$ $ \xck $ will
stand for
$<$ or $=,$ conversely, given $ \xck,$ $<$ will stand for $ \xck,$ but
not $=,$ similarly
for $ \xeb $ etc.
\index{Definition Log-Base}

\ed

\bd

$\hspace{0.01em}$

% (+++*** Orig. No.:  Definition Log-Base )

\label{Definition Log-Base}

We work here in a classical propositional language $ \xdl,$ a theory $T$
will be an
arbitrary set of formulas. Formulas will often be named $ \xbf,$ $ \xbq
,$ etc., theories
$T,$ $S,$ etc.

$v( \xdl )$ will be the set of propositional variables of $ \xdl.$

$M_{ \xdl }$ will be the set of (classical) models of $ \xdl,$ $M(T)$ or
$M_{T}$
is the set of models of $T,$ likewise $M( \xbf )$ for a formula $ \xbf.$

$ \xdD_{ \xdl }:=\{M(T):$ $T$ a theory in $ \xdl \},$ the set of definable
model sets.

Note that, in classical propositional logic, $ \xCQ,M_{ \xdl } \xbe
\xdD_{ \xdl },$ $ \xdD_{ \xdl }$ contains
singletons, is closed under arbitrary intersections and finite unions.

An operation $f: \xdy \xcp \xdp (M_{ \xdl })$ for $ \xdy \xcc \xdp (M_{
\xdl })$ is called definability
preserving, $ \xCf (dp)$ or $( \xbm dp)$ in short, iff for all $X \xbe
\xdD_{ \xdl } \xcs \xdy $ $f(X) \xbe \xdD_{ \xdl }.$

We will also use $( \xbm dp)$ for binary functions $f: \xdy \xCK \xdy \xcp
\xdp (M_{ \xdl })$ - as needed
for theory revision - with the obvious meaning.

$ \xcl $ will be classical derivability, and

$ \ol{T}:=\{ \xbf:T \xcl \xbf \},$ the closure of $T$ under $ \xcl.$

$Con(.)$ will stand for classical consistency, so $Con( \xbf )$ will mean
that
$ \xbf $ is clasical consistent, likewise for $Con(T).$ $Con(T,T' )$ will
stand for
$Con(T \xcv T' ),$ etc.

Given a consequence relation $ \xcn,$ we define

$ \ol{ \ol{T} }:=\{ \xbf:T \xcn \xbf \}.$

(There is no fear of confusion with $ \ol{T},$ as it just is not useful to
close
twice under classical logic.)

$T \xco T':=\{ \xbf \xco \xbf ': \xbf \xbe T, \xbf ' \xbe T' \}.$

If $X \xcc M_{ \xdl },$ then $Th(X):=\{ \xbf:X \xcm \xbf \},$ likewise
for $Th(m),$ $m \xbe M_{ \xdl }.$
\index{Definition Log-Cond}

\ed

\bd

$\hspace{0.01em}$

% (+++*** Orig. No.:  Definition Log-Cond )

\label{Definition Log-Cond}

We introduce here formally a list of properties of set functions on the
algebraic side, and their corresponding logical rules on the other side.

Recall that $ \ol{T}:=\{ \xbf:T \xcl \xbf \},$ $ \ol{ \ol{T} }:=\{ \xbf
:T \xcn \xbf \},$
where $ \xcl $ is classical consequence, and $ \xcn $ any other
consequence.

We show, wherever adequate, in parallel the formula version
in the left column, the theory version
in the middle column, and the semantical or algebraic
counterpart in the
right column. The algebraic counterpart gives conditions for a
function $f:\xdy\xcp\xdp (U)$, where $U$ is some set, and
$\xdy\xcc\xdp (U)$.

When the formula version is not commonly used, we omit it,
as we normally work only with the theory version.

Intuitively, $A$ and $B$ in the right hand side column stand for
$M(\xbf)$ for some formula $\xbf$, whereas $X$, $Y$ stand for
$M(T)$ for some theory $T$.

{\footnotesize

\begin{tabular}{|c|c|c|}

\hline

\multicolumn{3}{|c|}{Basics} \xEP

\hline

$(AND)$
\xEH
$(AND)$
\xEH
Closure under
\xEP

$ \xbf \xcn \xbq,  \xbf \xcn \xbq '   \xch $
\xEH
$ T \xcn \xbq, T \xcn \xbq '   \xch $
\xEH
finite
\xEP

$ \xbf \xcn \xbq \xcu \xbq ' $
\xEH
$ T \xcn \xbq \xcu \xbq ' $
\xEH
intersection
\xEP

\hline

$(OR)$ \xEH $(OR)$ \xEH $( \xbm OR)$ \xEP

$ \xbf \xcn \xbq,  \xbf ' \xcn \xbq   \xch $ \xEH
$ \ol{\ol{T}} \xcs \ol{\ol{T'}} \xcc \ol{\ol{T \xco T'}} $ \xEH
$f(X \xcv Y) \xcc f(X) \xcv f(Y)$
\xEP

$ \xbf \xco \xbf ' \xcn \xbq $ \xEH
\xEH
\xEP

\hline

$(wOR)$
\xEH
$(wOR)$
\xEH
$( \xbm wOR)$
\xEP

$ \xbf \xcn \xbq,$ $ \xbf ' \xcl \xbq $ $ \xch $
\xEH
$ \ol{ \ol{T} } \xcs \ol{T' }$ $ \xcc $ $ \ol{ \ol{T \xco T' } }$
\xEH
$f(X \xcv Y) \xcc f(X) \xcv Y$
\xEP

$ \xbf \xco \xbf ' \xcn \xbq $
\xEH
\xEH
\xEP

\hline

$(disjOR)$
\xEH
$(disjOR)$
\xEH
$( \xbm disjOR)$
\xEP

$ \xbf \xcl \xCN \xbf ',$ $ \xbf \xcn \xbq,$
\xEH
$\xCN Con(T \xcv T') \xch$
\xEH
$X \xcs Y= \xCQ $ $ \xch $
\xEP

$ \xbf ' \xcn \xbq $ $ \xch $ $ \xbf \xco \xbf ' \xcn \xbq $
\xEH
$ \ol{ \ol{T} } \xcs \ol{ \ol{T' } } \xcc \ol{ \ol{T \xco T' } }$
\xEH
$f(X \xcv Y) \xcc f(X) \xcv f(Y)$
\xEP

\hline

$(LLE)$
\xEH
$(LLE)$
\xEH
\xEP

Left Logical Equivalence
\xEH
\xEH
\xEP

$ \xcl \xbf \xcr \xbf ',  \xbf \xcn \xbq   \xch $
\xEH
$ \ol{T}= \ol{T' }  \xch   \ol{\ol{T}} = \ol{\ol{T'}}$
\xEH
trivially true
\xEP

$ \xbf ' \xcn \xbq $ \xEH \xEH \xEP

\hline

$(RW)$ Right Weakening
\xEH
$(RW)$
\xEH
upward closure
\xEP

$ \xbf \xcn \xbq,  \xcl \xbq \xcp \xbq '   \xch $
\xEH
$ T \xcn \xbq,  \xcl \xbq \xcp \xbq '   \xch $
\xEH
\xEP

$ \xbf \xcn \xbq ' $
\xEH
$T \xcn \xbq ' $
\xEH
\xEP

\hline

$(CCL)$ Classical Closure \xEH $(CCL)$ \xEH \xEP

\xEH
$ \ol{ \ol{T} }$ is classically
\xEH
trivially true
\xEP

\xEH closed \xEH \xEP

\hline

$(SC)$ Supraclassicality \xEH $(SC)$ \xEH $( \xbm \xcc )$ \xEP

$ \xbf \xcl \xbq $ $ \xch $ $ \xbf \xcn \xbq $ \xEH $ \ol{T} \xcc \ol{
\ol{T} }$ \xEH $f(X) \xcc X$ \xEP

\cline{1-1}

$(REF)$ Reflexivity \xEH \xEH \xEP
$ \xbD,\xba \xcn \xba $ \xEH \xEH \xEP

\hline

$(CP)$ \xEH $(CP)$ \xEH $( \xbm \xCQ )$ \xEP

Consistency Preservation \xEH \xEH \xEP

$ \xbf \xcn \xcT $ $ \xch $ $ \xbf \xcl \xcT $ \xEH $T \xcn \xcT $ $ \xch
$ $T \xcl \xcT $ \xEH $f(X)= \xCQ $ $ \xch $ $X= \xCQ $ \xEP

\hline

\xEH
\xEH $( \xbm \xCQ fin)$
\xEP

\xEH
\xEH $X \xEd \xCQ $ $ \xch $ $f(X) \xEd \xCQ $
\xEP

\xEH
\xEH for finite $X$
\xEP

\hline

\xEH $(PR)$ \xEH $( \xbm PR)$ \xEP

$ \ol{ \ol{ \xbf \xcu \xbf ' } }$ $ \xcc $ $ \ol{ \ol{ \ol{ \xbf } } \xcv
\{ \xbf ' \}}$ \xEH
$ \ol{ \ol{T \xcv T' } }$ $ \xcc $ $ \ol{ \ol{ \ol{T} } \xcv T' }$ \xEH
$X \xcc Y$ $ \xch $
\xEP

\xEH \xEH $f(Y) \xcs X \xcc f(X)$
\xEP

\cline{3-3}

\xEH
\xEH
$(\xbm PR ')$
\xEP

\xEH
\xEH
$f(X) \xcs Y \xcc f(X \xcs Y)$
\xEP

\hline

$(CUT)$ \xEH $(CUT)$ \xEH $ (\xbm CUT) $ \xEP
$ \xbD \xcn \xba; \xbD, \xba \xcn \xbb \xch $ \xEH
$T \xcc \ol{T' } \xcc \ol{ \ol{T} }  \xch $ \xEH
$f(X) \xcc Y \xcc X  \xch $ \xEP
$ \xbD \xcn \xbb $ \xEH
$ \ol{ \ol{T'} } \xcc \ol{ \ol{T} }$ \xEH
$f(X) \xcc f(Y)$
\xEP

\hline

% \end{tabular*}
\end{tabular}

}

{\footnotesize

% \begin{tabular*}{15.5cm}{|c@{\extracolsep\fill}|c|c|}
% \begin{tabular*}{15.23cm}{|c|c|c|}
\begin{tabular}{|c|c|c|}

\hline

\multicolumn{3}{|c|}{Cumulativity} \xEP

\hline

$(CM)$ Cautious Monotony \xEH $(CM)$ \xEH $ (\xbm CM) $ \xEP

$ \xbf \xcn \xbq,  \xbf \xcn \xbq '   \xch $ \xEH
$T \xcc \ol{T' } \xcc \ol{ \ol{T} }  \xch $ \xEH
$f(X) \xcc Y \xcc X  \xch $
\xEP

$ \xbf \xcu \xbq \xcn \xbq ' $ \xEH
$ \ol{ \ol{T} } \xcc \ol{ \ol{T' } }$ \xEH
$f(Y) \xcc f(X)$
\xEP

\cline{1-1}

\cline{3-3}

or $(ResM)$ Restricted Monotony \xEH \xEH $(\xbm ResM)$ \xEP
$ \xbD \xcn \xba, \xbb \xch \xbD,\xba \xcn \xbb $ \xEH \xEH
$ f(X) \xcc A \xcs B \xch f(X \xcs A) \xcc B $ \xEP

\hline

$(CUM)$ Cumulativity \xEH $(CUM)$ \xEH $( \xbm CUM)$ \xEP

$ \xbf \xcn \xbq   \xch $ \xEH
$T \xcc \ol{T' } \xcc \ol{ \ol{T} }  \xch $ \xEH
$f(X) \xcc Y \xcc X  \xch $
\xEP

$( \xbf \xcn \xbq '   \xcj   \xbf \xcu \xbq \xcn \xbq ' )$ \xEH
$ \ol{ \ol{T} }= \ol{ \ol{T' } }$ \xEH
$f(Y)=f(X)$ \xEP

\hline

\xEH
$ (\xcc \xcd) $
\xEH
$ (\xbm \xcc \xcd) $
\xEP
\xEH
$T \xcc \ol{\ol{T'}}, T' \xcc \ol{\ol{T}} \xch $
\xEH
$ f(X) \xcc Y, f(Y) \xcc X \xch $
\xEP
\xEH
$ \ol{\ol{T'}} = \ol{\ol{T}}$
\xEH
$ f(X)=f(Y) $
\xEP

\hline

\multicolumn{3}{|c|}{Rationality} \xEP

\hline

$(RatM)$ Rational Monotony \xEH $(RatM)$ \xEH $( \xbm RatM)$ \xEP

$ \xbf \xcn \xbq,  \xbf \xcN \xCN \xbq '   \xch $ \xEH
$Con(T \xcv \ol{\ol{T'}})$, $T \xcl T'$ $ \xch $ \xEH
$X \xcc Y, X \xcs f(Y) \xEd \xCQ   \xch $
\xEP

$ \xbf \xcu \xbq ' \xcn \xbq $ \xEH
$ \ol{\ol{T}} \xcd \ol{\ol{\ol{T'}} \xcv T} $ \xEH
$f(X) \xcc f(Y) \xcs X$ \xEP

\hline

\xEH $(RatM=)$ \xEH $( \xbm =)$ \xEP

\xEH
$Con(T \xcv \ol{\ol{T'}})$, $T \xcl T'$ $ \xch $ \xEH
$X \xcc Y, X \xcs f(Y) \xEd \xCQ   \xch $
\xEP

\xEH
$ \ol{\ol{T}} = \ol{\ol{\ol{T'}} \xcv T} $ \xEH
$f(X) = f(Y) \xcs X$ \xEP

\hline

\xEH
$(Log=' )$
\xEH $( \xbm =' )$
\xEP

\xEH
$Con( \ol{ \ol{T' } } \xcv T)$ $ \xch $
\xEH $f(Y) \xcs X \xEd \xCQ $ $ \xch $
\xEP

\xEH
$ \ol{ \ol{T \xcv T' } }= \ol{ \ol{ \ol{T' } } \xcv T}$
\xEH $f(Y \xcs X)=f(Y) \xcs X$
\xEP

\hline

\xEH
$(Log \xFO )$
\xEH $( \xbm \xFO )$
\xEP

\xEH
$ \ol{ \ol{T \xco T' } }$ is one of
\xEH $f(X \xcv Y)$ is one of
\xEP

\xEH
$\ol{\ol{T}},$ or $\ol{\ol{T'}},$ or $\ol{\ol{T}} \xcs \ol{\ol{T'}}$ (by (CCL))
\xEH $f(X),$ $f(Y)$ or $f(X) \xcv f(Y)$
\xEP

\hline

\xEH
$(Log \xcv )$
\xEH $( \xbm \xcv )$
\xEP

\xEH
$Con( \ol{ \ol{T' } } \xcv T),$ $ \xCN Con( \ol{ \ol{T' } }
\xcv \ol{ \ol{T} })$ $ \xch $
\xEH $f(Y) \xcs (X-f(X)) \xEd \xCQ $ $ \xch $
\xEP

\xEH
$ \xCN Con( \ol{ \ol{T \xco T' } } \xcv T' )$
\xEH $f(X \xcv Y) \xcs Y= \xCQ$
\xEP

\hline

\xEH
$(Log \xcv ' )$
\xEH $( \xbm \xcv ' )$
\xEP

\xEH
$Con( \ol{ \ol{T' } } \xcv T),$ $ \xCN Con( \ol{ \ol{T' }
} \xcv \ol{ \ol{T} })$ $ \xch $
\xEH $f(Y) \xcs (X-f(X)) \xEd \xCQ $ $ \xch $
\xEP

\xEH
$ \ol{ \ol{T \xco T' } }= \ol{ \ol{T} }$
\xEH $f(X \xcv Y)=f(X)$
\xEP

\hline

\xEH
\xEH $( \xbm \xbe )$
\xEP

\xEH
\xEH $a \xbe X-f(X)$ $ \xch $
\xEP

\xEH
\xEH $ \xcE b \xbe X.a \xce f(\{a,b\})$
\xEP

\hline

\end{tabular}

}

$(PR)$ is also called infinite conditionalization - we choose the name for
its central role for preferential structures $(PR)$ or $( \xbm PR).$

The system of rules $(AND)$ $(OR)$ $(LLE)$ $(RW)$ $(SC)$ $(CP)$ $(CM)$ $(CUM)$
is also called system $P$ (for preferential), adding $(RatM)$ gives the system
$R$ (for rationality or rankedness).

Roughly: Smooth preferential structures generate logics satisfying system
$P$, ranked structures logics satisfying system $R$.

A logic satisfying $(REF)$, $(ResM)$, and $(CUT)$ is called a consequence
relation.

$(LLE)$ and$(CCL)$ will hold automatically, whenever we work with model sets.

$(AND)$ is obviously closely related to filters, and corresponds to closure
under finite intersections. $(RW)$ corresponds to upward closure of filters.

More precisely, validity of both depend on the definition, and the
direction we consider.

Given $f$ and $(\xbm \xcc )$, $f(X)\xcc X$ generates a pricipal filter:
$\{X'\xcc X:f(X)\xcc X'\}$, with
the definition: If $X=M(T)$, then $T\xcn \xbf$  iff $f(X)\xcc M(\xbf )$.
Validity of $(AND)$ and
$(RW)$ are then trivial.

Conversely, we can define for $X=M(T)$

$\xdx:=\{X'\xcc X: \xcE \xbf (X'=X\xcs M(\xbf )$ and $T\xcn \xbf )\}$.

$(AND)$ then makes $\xdx$  closed under
finite intersections, $(RW)$ makes $\xdx$  upward
closed. This is in the infinite case usually not yet a filter, as not all
subsets of $X$ need to be definable this way.
In this case, we complete $\xdx$  by
adding all $X''$ such that there is $X'\xcc X''\xcc X$, $X'\xbe\xdx$.

Alternatively, we can define

$\xdx:=\{X'\xcc X: \xcS\{X \xcs M(\xbf ): T\xcn \xbf \} \xcc X' \}$.

$(SC)$ corresponds to the choice of a subset.

$(CP)$ is somewhat delicate, as it presupposes that the chosen model set is
non-empty. This might fail in the presence of ever better choices, without
ideal ones; the problem is addressed by the limit versions.

$(PR)$ is an infinitary version of one half of the deduction theorem: Let $T$
stand for $\xbf$, $T'$ for $\xbq$, and $\xbf \xcu \xbq \xcn \xbs$,
so $\xbf \xcn \xbq \xcp \xbs$, but $(\xbq \xcp \xbs )\xcu \xbq \xcl \xbs$.

$(CUM)$ (whose most interesting half in our context is $(CM)$) may best be seen
as
normal use of lemmas: We have worked hard and found some lemmas. Now
we can take a rest, and come back again with our new lemmas. Adding them to the
axioms will neither add new theorems, nor prevent old ones to hold.

\index{Definition Pref-Str}

\ed

\bd

$\hspace{0.01em}$

% (+++*** Orig. No.:  Definition Pref-Str )

\label{Definition Pref-Str}

Fix $U \xEd \xCQ,$ and consider arbitrary $X.$
Note that this $X$ has not necessarily anything to do with $U,$ or $ \xdu
$ below.
Thus, the functions $ \xbm_{ \xdm }$ below are in principle functions from
$V$ to $V$ - where $V$
is the set theoretical universe we work in.

(A) Preferential models or structures.

(1) The version without copies:

A pair $ \xdm:=<U, \xeb >$ with $U$ an arbitrary set, and $ \xeb $ an
arbitrary binary relation
is called a preferential model or structure.

(2) The version with copies:

A pair $ \xdm:=< \xdu, \xeb >$ with $ \xdu $ an arbitrary set of pairs,
and $ \xeb $ an arbitrary binary
relation is called a preferential model or structure.

If $<x,i> \xbe \xdu,$ then $x$ is intended to be an element of $U,$ and
$i$ the index of the
copy.

We sometimes also need copies of the relation $ \xeb,$ we will then
replace $ \xeb $
by one or several arrows $ \xba $ attacking non-minimal elements, e.g. $x
\xeb y$ will
be written $ \xba:x \xcp y,$ $<x,i> \xeb <y,i>$ will be written $ \xba
:<x,i> \xcp <y,i>,$ and
finally we might have $< \xba,k>:x \xcp y$ and $< \xba,k>:<x,i> \xcp
<y,i>,$ etc.

(B) Minimal elements, the functions $ \xbm_{ \xdm }$

(1) The version without copies:

Let $ \xdm:=<U, \xeb >,$ and define

$ \xbm_{ \xdm }(X)$ $:=$ $\{x \xbe X:$ $x \xbe U$ $ \xcu $ $ \xCN \xcE x'
\xbe X \xcs U.x' \xeb x\}.$

$ \xbm_{ \xdm }(X)$ is called the set of minimal elements of $X$ (in $
\xdm ).$

(2) The version with copies:

Let $ \xdm:=< \xdu, \xeb >$ be as above. Define

$ \xbm_{ \xdm }(X)$ $:=$ $\{x \xbe X:$ $ \xcE <x,i> \xbe \xdu. \xCN \xcE
<x',i' > \xbe \xdu (x' \xbe X$ $ \xcu $ $<x',i' >' \xeb <x,i>)\}.$

Again, by abuse of language, we say that $ \xbm_{ \xdm }(X)$ is the set of
minimal elements
of $X$ in the structure. If the context is clear, we will also write just
$ \xbm.$

We sometimes say that $<x,i>$ ``kills'' or ``minimizes'' $<y,j>$ if
$<x,i> \xeb <y,j>.$ By abuse of language we also say a set $X$ kills or
minimizes a set
$Y$ if for all $<y,j> \xbe \xdu,$ $y \xbe Y$ there is $<x,i> \xbe \xdu,$
$x \xbe X$ s.t. $<x,i> \xeb <y,j>.$

$ \xdm $ is also called injective or 1-copy, iff there is always at most
one copy
$<x,i>$ for each $x.$ Note that the existence of copies corresponds to a
non-injective labelling function - as is often used in nonclassical
logic, e.g. modal logic.

We say that $ \xdm $ is transitive, irreflexive, etc., iff $ \xeb $ is.

Note that $ \xbm (X)$ might well be empty, even if $X$ is not.
\index{Definition Pref-Log}

\ed

\bd

$\hspace{0.01em}$

% (+++*** Orig. No.:  Definition Pref-Log )

\label{Definition Pref-Log}

We define the consequence relation of a preferential structure for a
given propositional language $ \xdl.$

(A)

(1) If $m$ is a classical model of a language $ \xdl,$ we say by abuse
of language

$<m,i> \xcm \xbf $ iff $m \xcm \xbf,$

and if $X$ is a set of such pairs, that

$X \xcm \xbf $ iff for all $<m,i> \xbe X$ $m \xcm \xbf.$

(2) If $ \xdm $ is a preferential structure, and $X$ is a set of $ \xdl
-$models for a
classical propositional language $ \xdl,$ or a set of pairs $<m,i>,$
where the $m$ are
such models, we call $ \xdm $ a classical preferential structure or model.

(B)

Validity in a preferential structure, or the semantical consequence
relation
defined by such a structure:

Let $ \xdm $ be as above.

We define:

$T \xcm_{ \xdm } \xbf $ iff $ \xbm_{ \xdm }(M(T)) \xcm \xbf,$ i.e. $
\xbm_{ \xdm }(M(T)) \xcc M( \xbf ).$

$ \xdm $ will be called definability preserving iff for all $X \xbe \xdD_{
\xdl }$ $ \xbm_{ \xdm }(X) \xbe \xdD_{ \xdl }.$

As $ \xbm_{ \xdm }$ is defined on $ \xdD_{ \xdl },$ but need by no means
always result in some new
definable set, this is (and reveals itself as a quite strong) additional
property.
\index{Definition Smooth}

\ed

\bd

$\hspace{0.01em}$

% (+++*** Orig. No.:  Definition Smooth )

\label{Definition Smooth}

Let $ \xdy \xcc \xdp (U).$ (In applications to logic, $ \xdy $ will be $
\xdD_{ \xdl }.)$

A preferential structure $ \xdm $ is called $ \xdy -$smooth iff in every
$X \xbe \xdy $ every element
$x \xbe X$ is either minimal in $X$ or above an element, which is minimal
in $X.$ More
precisely:

(1) The version without copies:

If $x \xbe X \xbe \xdy,$ then either $x \xbe \xbm (X)$ or there is $x'
\xbe \xbm (X).x' \xeb x.$

(2) The version with copies:

If $x \xbe X \xbe \xdy,$ and $<x,i> \xbe \xdu,$ then either there is no
$<x',i' > \xbe \xdu,$ $x' \xbe X,$
$<x',i' > \xeb <x,i>$ or there is $<x',i' > \xbe \xdu,$ $<x',i' > \xeb
<x,i>,$ $x' \xbe X,$ s.t. there is
no $<x'',i'' > \xbe \xdu,$ $x'' \xbe X,$ with $<x'',i'' > \xeb <x',i'
>.$

When considering the models of a language $ \xdl,$ $ \xdm $ will be
called smooth iff
it is $ \xdD_{ \xdl }-$smooth; $ \xdD_{ \xdl }$ is the default.

Obviously, the richer the set $ \xdy $ is, the stronger the condition $
\xdy -$smoothness
will be.

% karl-search= End Hier-ARepr-Def
\vspace{7mm}

% *************************************

\vspace{7mm}

\index{Fact Rank-Base}

\ed

\bfa

$\hspace{0.01em}$

% (+++*** Orig. No.:  Fact Rank-Base )

\label{Fact Rank-Base}

Let $ \xeb $ be an irreflexive, binary relation on $X,$ then the following
two conditions
are equivalent:

(1) There is $ \xbO $ and an irreflexive, total, binary relation $ \xeb '
$ on $ \xbO $ and a
function $f:X \xcp \xbO $ s.t. $x \xeb y$ $ \xcr $ $f(x) \xeb ' f(y)$ for
all $x,y \xbe X.$

(2) Let $x,y,z \xbe X$ and $x \xcT y$ wrt. $ \xeb $ (i.e. neither $x \xeb
y$ nor $y \xeb x),$ then $z \xeb x$ $ \xcp $ $z \xeb y$
and $x \xeb z$ $ \xcp $ $y \xeb z.$

$ \xcz $
\\[3ex]
\index{Definition Rank-Rel}

\efa

\bd

$\hspace{0.01em}$

% (+++*** Orig. No.:  Definition Rank-Rel )

\label{Definition Rank-Rel}

We call an irreflexive, binary relation $ \xeb $ on $X,$ which satisfies
(1)
(equivalently (2)) of Fact \ref{Fact Rank-Base}, ranked.
By abuse of language, we also call a preferential structure $<X, \xeb >$
ranked, iff
$ \xeb $ is.

% {\LARGE karl-search= Start Hier-Def }

\label{Section Hier-Def}
\index{Section Hier-Def}

\ed

\bd

$\hspace{0.01em}$

% (+++*** Orig. No.:  Definition 2.1 )

\label{Definition 2.1}

We have the usual framework of preferential structures, i.e.
either a set with a possibly non-injective labelling function, or,
equivalently,
a set of possible worlds with copies. The relation of the
preferential structure will be fixed, and will not depend on the point $m$
from
where we look at it.

Next, we have a set $ \xdA,$ and a finite, disjoint cover $A_{i}:i<n$ of
$ \xda,$
with a relation ``of quality'' $<,$ $ \xda $ will denote the $A_{i}$ (and
thus $ \xdA ),$ and $<,$
i.e. $ \xda =<\{A_{i}:i \xbe I\},<>.$

By Fact \ref{Fact 6.21}, we may assume that all $A_{i}$ are described
by a formula.

Finally, we have $ \xdB \xcc \xdA,$ the subset of ``good'' elements of $
\xdA $ - which
we also assume to be described by a formula.

In addition, we have a binary relation of accessibility, $R,$ which we
assume
transitive - modal operators will be defined relative to $R.$ $R$
determines which
part of the preferential structure is visible.

Let $R(s):=\{t:sRt\}.$

\ed

\bd

$\hspace{0.01em}$

% (+++*** Orig. No.:  Definition 2.2 )

\label{Definition 2.2}

We repeat here from the introduction, and assume $A_{i}=M( \xba_{i}),$
$B=M( \xbb ),$ and
$ \xbm $ expresses the minimality of the preferential structure.

$t \xcm \xba_{i}> \xbb: \xcj $ $ \xbm (A_{i}) \xcs R(t) \xcc B,$

we will also abuse notation and just write

$t \xcm A_{i}>B$ in this case.

We then define:

$t \xcm \xdc $ iff at the smallest $i$ s.t. $ \xbm (A_{i}) \xcs R(t) \xEd
\xCQ,$ $ \xbm (A_{i}) \xcs R(t) \xcc \xdB $ holds.

\ed

This motivates the following:

\bd

$\hspace{0.01em}$

% (+++*** Orig. No.:  Definition 2.7 )

\label{Definition 2.7}

Let $ \xdA $ be a fixed set, and $ \xda $ a finite, totally ordered (by
$<)$ disjoint cover
by non-empty subsets of $ \xdA.$

For $x \xbe \xdA,$ let $rg(x)$ the unique $A \xbe \xda $ such that $x
\xbe A,$ so $rg(x)<rg(y)$ is
defined in the natural way.

A preferential structure $< \xdx, \xeb >$ $( \xdx $ a set of pairs
$<x,i>)$ is called $ \xda -$ranked
iff for all $x,x' $ $rg(x)<rg(x' )$ implies $<x,i> \xeb <x',i' >$ for all
$<x,i>,<x',i' > \xbe \xdx.$

\ed

Note that automatically for $X \xcc \xdA,$ $ \xbm (X) \xcc A_{j}$ when
$j$ is the smallest $i$ s.t.
$X \xcs A_{i} \xEd \xCQ.$

The idea is now to make the $A_{i}$ the layers, and ``trigger'' the first
layer
$A_{j}$ s.t. $ \xbm (A_{j}) \xcs R(x) \xEd \xCQ,$ and check whether $
\xbm (A_{j}) \xcs R(x) \xcc B_{j}.$ A suitable
ranked structure will automatically find this $A_{j}.$

More definitions and results for such $ \xda $ and $ \xdc $ will be found
in Section \ref{Section Hier-CondRepr}.

% karl-search= End Hier-Def
\vspace{7mm}

% *************************************

\vspace{7mm}

%  6.2  Discussion
%  6.2  Discussion
% %
% ================
\subsection{
Discussion
}

% {\LARGE karl-search= Start Hier-ARepr-Disc }

\label{Section Hier-ARepr-Disc}
\index{Section Hier-ARepr-Disc}

The not necessarily smooth and the smooth case will be treated
differently.

Strangely, the smooth case is simpler, as an added new layer in the proof
settles it. Yet, this is not surprising when looking closer, as minimal
elements
never have
higher rank, and we know from $( \xbm CUM)$ that minimizing by minimal
elements
suffices. All we have to add that any element in the minimal layer
minimizes
any element higher up.

In the simple, not necessarily smooth, case, we have to go deeper into the
original proof to obtain the result.

The following idea, inspired by the treatment of the smooth case, will not
work: Instead of minimizing by arbitrary elements, minimize only by
elements of
minimal rank, as the following example shows. If it worked, we might add
just
another layer to the original proof without $( \xbm \xda ),$
(see Definition \ref{Definition 6.6}), as in the smooth case.

\be

$\hspace{0.01em}$

% (+++*** Orig. No.:  Example 6.1 )

\label{Example 6.1}

Consider the base set $\{a,b,c\},$ $ \xbm (\{a,b,c\})=\{b\},$ $ \xbm
(\{a,b\})=\{a,b\},$
$ \xbm (\{a,c\})= \xCQ,$ $ \xbm (\{b,c\})=\{b\},$ $ \xda $ defined by
$\{a,b\}<\{c\}.$

Obviously, $( \xbm \xda )$ is satisfied. $ \xbm $ can be represented by
the (not transitive!)
relation $a \xeb c \xeb a,$ $b \xeb c,$ which is $ \xda -$ranked.

But trying to minimize a in $\{a,b,c\}$ in the minimal layer will lead to
$b \xeb a,$
and thus $a \xce \xbm (\{a,b\}),$ which is wrong.

$ \xcz $
\\[3ex]

\ee

Both proofs are self contained, the proof of the smooth case is stronger
than
the one published in  \cite{Sch04}, as it does not need closure
under finite
intersections. The proofs of the general and transitive general case are
adaptations of earlier proofs by the second author, but the basic ideas
are not new, and were published before - see e.g.  \cite{Sch04}, or
 \cite{Sch92}. The proofs for the smooth case are slightly
stronger than those published before (see again e.g.  \cite{Sch04}),
as we
work without closure under finite intersection. The rest is almost
verbatim the
same, and we only add a supplementary layer
in the end (Fact \ref{Fact Smooth-to-A}), which
will make the construction $ \xda -$ranked.

In the following, we will assume the partition $ \xda $ to be given. We
could also
construct it from the properties of $ \xbm,$ but this would need stronger
closure
properties of the domain. The construction of $ \xda $ is more difficult
than the
construction of the ranking in fully ranked structures, as $x \xbe \xbm
(X),$ $y \xbe X- \xbm (X)$
will guarantee only $rg(x) \xck rg(y),$ and not $rg(x)<rg(y),$ as is the
case in the
latter situation. This corresponds to the separate treatment of the $ \xba
$ and
other formulas in the logical version, discussed
in Section \ref{Section Hier-ARepr-Logic}.

% karl-search= End Hier-ARepr-Disc
\vspace{7mm}

% *************************************

\vspace{7mm}

%  6.3  &a-ranked general and transitive structures
%  6.3  &a-ranked general and transitive structures
% %
% =================================================
\subsection{
$\xda$-ranked general and transitive structures
}

% {\LARGE karl-search= Start Hier-ARepr-General }

\label{Section Hier-ARepr-General}
\index{Section Hier-ARepr-General}

We will show here the following representation results:

Let $ \xda $ be given.

An operation $ \xbm: \xdy \xcp \xdp (Z)$ is representable by an $ \xda
-$ranked preferential
structure iff $ \xbm $ satisfies $( \xbm \xcc ),$ $( \xbm PR),$ $( \xbm
\xda )$ (Proposition \ref{Proposition 6.3}), and,
moreover, the structure can be chosen transitive (Proposition \ref{Proposition
6.5}).

Note that we carefully avoid any unnecessary assumptions about the domain
$ \xdy \xcc \xdp (Z)$ of the function $ \xbm.$

\bd

$\hspace{0.01em}$

% (+++*** Orig. No.:  Definition 6.6 )

\label{Definition 6.6}

We define a new condition:

Let $ \xda $ be given as defined in Definition \ref{Definition 2.6}.

$( \xbm \xda )$ If $X \xbe \xdy,$ $A,A' \xbe \xda,$ $A<A',$ $X \xcs A
\xEd \xCQ,$ $X \xcs A' \xEd \xCQ $ then $ \xbm (X) \xcs A' = \xCQ.$

\ed

This new condition will be central for the modified representation.
%  6.3.1  The basic, not necessarily transitive, case
%  6.3.1  The basic, not necessarily transitive, case
% %
% ===================================================
\subsubsection{
The basic, not necessarily transitive, case
}

% {\LARGE karl-search= Start Hier-ARepr-General-Intrans }

\label{Section Hier-ARepr-General-Intrans}
\index{Section Hier-ARepr-General-Intrans}
\index{Definition Y-Pi-x}

\bd

$\hspace{0.01em}$

% (+++*** Orig. No.:  Definition Y-Pi-x )

\label{Definition Y-Pi-x}

For $x \xbe Z,$ let $ \xdy_{x}:=\{Y \xbe \xdy $: $x \xbe Y- \xbm (Y)\},$
$ \xbP_{x}:= \xbP \xdy_{x}.$

\ed

Note that $ \xCQ \xce \xdy_{x}$, $ \xbP_{x} \xEd \xCQ,$ and that $
\xbP_{x}=\{ \xCQ \}$ iff $ \xdy_{x}= \xCQ.$
\index{Claim Mu-f}

\bc

$\hspace{0.01em}$

% (+++*** Orig. No.:  Claim Mu-f )

\label{Claim Mu-f}

Let $ \xbm: \xdy \xcp \xdp (Z)$ satisfy $( \xbm \xcc )$ and $( \xbm PR),$
and let $U \xbe \xdy.$
Then $x \xbe \xbm (U)$ $ \xcr $ $x \xbe U$ $ \xcu $ $ \xcE f \xbe
\xbP_{x}.ran(f) \xcs U= \xCQ.$
\index{Claim Mu-f Proof}

\ec

\subparagraph{
Proof:
}

$\hspace{0.01em}$

% (+++*** Orig.:  Proof: )

\label{Section Proof:}

Case 1: $ \xdy_{x}= \xCQ,$ thus $ \xbP_{x}=\{ \xCQ \}.$
`` $ \xcp $ '': Take $f:= \xCQ.$
`` $ \xcq $ '': $x \xbe U \xbe \xdy,$ $ \xdy_{x}= \xCQ $ $ \xcp $ $x \xbe
\xbm (U)$ by definition of $ \xdy_{x}.$

Case 2: $ \xdy_{x} \xEd \xCQ.$
`` $ \xcp $ '': Let $x \xbe \xbm (U) \xcc U.$ It suffices to show $Y \xbe
\xdy_{x}$ $ \xcp $ $Y-U \xEd \xCQ.$ But if $Y \xcc U$ and
$Y \xbe \xdy_{x}$, then $x \xbe Y- \xbm (Y),$ contradicting $( \xbm PR).$
`` $ \xcq $ '': If $x \xbe U- \xbm (U),$ then $U \xbe \xdy_{x}$, so $
\xcA f \xbe \xbP_{x}.ran(f) \xcs U \xEd \xCQ.$
$ \xcz $
\\[3ex]
\index{Construction Pref-Base}

\bcs

$\hspace{0.01em}$

% (+++*** Orig. No.:  Construction Pref-Base )

\label{Construction Pref-Base}

Let $ \xdx:=\{<x,f>:x \xbe Z$ $ \xcu $ $f \xbe \xbP_{x}\},$ and $<x',f'
> \xeb <x,f>$ $: \xcr $ $x' \xbe ran(f).$ Let $ \xdz:=< \xdx, \xeb >.$

\ecs

\bc

$\hspace{0.01em}$

% (+++*** Orig. No.:  Claim 6.1 )

\label{Claim 6.1}

Let $ \xbm: \xdy \xcp \xdp (Z)$ satisfy $( \xbm \xcc )$ and $( \xbm PR),$
and let $U \xbe \xdy.$
Then $x \xbe \xbm (U)$ $ \xcr $ $x \xbe U$ $ \xcu $ $ \xcE f \xbe
\xbP_{x}.ran(f) \xcs U= \xCQ.$

\ec

\subparagraph{
Proof:
}

$\hspace{0.01em}$

% (+++*** Orig.:  Proof: )

\label{Section Proof:}

Case 1: $ \xdy_{x}= \xCQ,$ thus $ \xbP_{x}=\{ \xCQ \}.$
`` $ \xcp $ '': Take $f:= \xCQ.$
`` $ \xcq $ '': $x \xbe U \xbe \xdy,$ $ \xdy_{x}= \xCQ $ $ \xcp $ $x \xbe
\xbm (U)$ by definition of $ \xdy_{x}.$

Case 2: $ \xdy_{x} \xEd \xCQ.$
`` $ \xcp $ '': Let $x \xbe \xbm (U) \xcc U.$ It suffices to show $Y \xbe
\xdy_{x}$ $ \xcp $ $Y-U \xEd \xCQ.$ But if $Y \xcc U$ and
$Y \xbe \xdy_{x}$, then $x \xbe Y- \xbm (Y),$ contradicting $( \xbm PR).$
`` $ \xcq $ '': If $x \xbe U- \xbm (U),$ then $U \xbe \xdy_{x}$, so $
\xcA f \xbe \xbP_{x}.ran(f) \xcs U \xEd \xCQ.$
$ \xcz $
\\[3ex]

\bco

$\hspace{0.01em}$

% (+++*** Orig. No.:  Corollary 6.2 )

\label{Corollary 6.2}

Let $ \xbm: \xdy \xcp \xdp (Z)$ satisfy $( \xbm \xcc ),$ $( \xbm PR),$ $(
\xbm \xda ),$ and let $U \xbe \xdy.$

If $x \xbe U$ and $ \xcE x' \xbe U.rg(x' )<rg(x),$ then $ \xcA f \xbe
\xbP_{x}.ran(f) \xcs U \xEd \xCQ.$

\eco

\subparagraph{
Proof:
}

$\hspace{0.01em}$

% (+++*** Orig.:  Proof: )

\label{Section Proof:}

By $( \xbm \xda )$ $x \xce \xbm (U),$ thus by Claim \ref{Claim 6.1} $
\xcA f \xbe \xbP_{x}.ran(f) \xcs U \xEd \xCQ.$ $ \xcz $
\\[3ex]

\bp

$\hspace{0.01em}$

% (+++*** Orig. No.:  Proposition 6.3 )

\label{Proposition 6.3}

Let $ \xda $ be given.

An operation $ \xbm: \xdy \xcp \xdp (Z)$ is representable by an $ \xda
-$ranked preferential
structure iff $ \xbm $ satisfies $( \xbm \xcc ),$ $( \xbm PR),$ $( \xbm
\xda ).$

\ep

\subparagraph{
Proof:
}

$\hspace{0.01em}$

% (+++*** Orig.:  Proof: )

\label{Section Proof:}

One direction is trivial. The central argument is: If $a \xeb b$ in $X,$
and $X \xcc Y,$
then $a \xeb b$ in $Y,$ too.

We turn to the other direction. The preferential structure is defined in
Construction \ref{Construction 6.1}, Claim \ref{Claim 6.4} shows
representation.

\bcs

$\hspace{0.01em}$

% (+++*** Orig. No.:  Construction 6.1 )

\label{Construction 6.1}

Let $ \xdx:=\{<x,f>:x \xbe Z$ $ \xcu $ $f \xbe \xbP_{x}\},$ and $<x',f'
> \xeb <x,f>$ $: \xcr $ $x' \xbe ran(f)$ or $rg(x' )<rg(x).$

Note that, as $ \xda $ is given, we also know $rg(x).$

Let $ \xdz:=< \xdx, \xeb >.$

\ecs

Obviously, $ \xdz $ is $ \xda -$ranked.

\bc

$\hspace{0.01em}$

% (+++*** Orig. No.:  Claim 6.4 )

\label{Claim 6.4}

For $U \xbe \xdy,$ $ \xbm (U)= \xbm_{ \xdz }(U).$

\ec

\subparagraph{
Proof:
}

$\hspace{0.01em}$

% (+++*** Orig.:  Proof: )

\label{Section Proof:}

By Claim \ref{Claim 6.1}, it suffices to show that for all $U \xbe
\xdy $
$x \xbe \xbm_{ \xdz }(U)$ $ \xcr $ $x \xbe U$ and $ \xcE f \xbe
\xbP_{x}.ran(f) \xcs U= \xCQ.$ So let $U \xbe \xdy.$

`` $ \xcp $ '': If $x \xbe \xbm_{ \xdz }(U),$ then there is $<x,f>$
minimal in $ \xdx \xex U$ - where
$ \xdx \xex U:=\{<x,i> \xbe \xdx:x \xbe U\}),$ so $x \xbe U,$ and there
is no $<x',f' > \xeb <x,f>,$ $x' \xbe U,$ so by $ \xbP_{x' } \xEd \xCQ $
there is no $x' \xbe ran(f),$ $x' \xbe U,$ but then
$ran(f) \xcs U= \xCQ.$

`` $ \xcq $ '': If $x \xbe U,$ and there is $f \xbe \xbP_{x}$, $ran(f)
\xcs U= \xCQ,$
then by Corollary \ref{Corollary 6.2}, there is no $x' \xbe U,$ $rg(x'
)<rg(x),$ so $<x,f>$ is
minimal in $ \xdx \xex U.$

$ \xcz $ (Claim \ref{Claim 6.4} and Proposition \ref{Proposition 6.3})
\\[3ex]

% karl-search= End Hier-ARepr-General-Intrans
\vspace{7mm}

% *************************************

\vspace{7mm}

%  6.3.2  The transitive case
%  6.3.2  The transitive case
% %
% ===========================
\subsubsection{
The transitive case
}

% {\LARGE karl-search= Start Hier-ARepr-General-Trans }

\label{Section Hier-ARepr-General-Trans}
\index{Section Hier-ARepr-General-Trans}

\bp

$\hspace{0.01em}$

% (+++*** Orig. No.:  Proposition 6.5 )

\label{Proposition 6.5}

Let $ \xda $ be given.

An operation $ \xbm: \xdy \xcp \xdp (Z)$ is representable by an $ \xda
-$ranked transitive
preferential structure iff $ \xbm $ satisfies $( \xbm \xcc ),$ $( \xbm
PR),$ $( \xbm \xda ).$

\ep

\bcs

$\hspace{0.01em}$

% (+++*** Orig. No.:  Construction 6.2 )

\label{Construction 6.2}

\index{$T_x$}
(1) For $x \xbe Z,$ let $T_{x}$ be the set of trees $t_{x}$ s.t.

(a) all nodes are elements of $Z,$

(b) the root of $t_{x}$ is $x,$

(c) $height(t_{x}) \xck \xbo,$

(d) if $y$ is an element in $t_{x}$, then there is $f \xbe \xbP_{y}:=
\xbP \{Y \xbe \xdy $: $y \xbe Y- \xbm (Y)\}$
s.t. the set of children of $y$ is $ran(f) \xcv \{y' \xbe Z:rg(y'
)<rg(y)\}.$

(2) For $x,y \xbe Z,$ $t_{x} \xbe T_{x}$, $t_{y} \xbe T_{y}$, set $t_{x}
\xem t_{y}$ iff $y$ is a (direct) child
of the root $x$ in $t_{x}$, and $t_{y}$ is the subtree of $t_{x}$
beginning at $y.$

(3) Let $ \xdz $ $:=$ $<$ $\{<x,t_{x}>:$ $x \xbe Z,$ $t_{x} \xbe T_{x}\}$
, $<x,t_{x}> \xee <y,t_{y}>$ iff $t_{x} \xem t_{y}$ $>.$

\ecs

\bfa

$\hspace{0.01em}$

% (+++*** Orig. No.:  Fact 6.6 )

\label{Fact 6.6}

(1) The construction ends at some $y$ iff $ \xdy_{y}= \xCQ $ and there is
no $y' $ s.t.
$rg(y' )<rg(y),$ consequently $T_{x}=\{x\}$ iff $ \xdy_{x}= \xCQ $ and
there are no $x' $ with
lesser rang. (We identify the tree of height 1 with its root.)

(2) We define a special tree $tc_{x}$ for all $x:$ For all nodes $y$ in
$tc_{x},$ the
successors are as follows:
if $ \xdy_{y} \xEd \xCQ,$ then $z$ is an successor iff $z=y$ or
$rg(z)<rg(y);$
if $ \xdy_{y}= \xCQ,$ then $z$ is an successor iff $rg(z)<rg(y).$
(In the first case, we make $f \xbe \xdy_{y}$ always choose $y$ itself.)
$tc_{x}$ is an element of $T_{x}.$ Thus, with (1), $T_{x} \xEd \xCQ $ for
any $x.$
Note: $tc_{x}=x$ iff $ \xdy_{x}= \xCQ $ and $x$ has minimal rang.

(3) If $f \xbe \xbP_{x}$, then the tree $tf_{x}$ with root $x$ and
otherwise
composed of the subtrees $tc_{y}$ for $y \xbe ran(f) \xcv \{y':rg(y'
)<rg(y)\}$
is an element of $T_{x}$.
(Level 0 of $tf_{x}$ has $x$ as element, the $t_{y}' s$ begin at level 1.)

(4) If $y$ is an element in $t_{x}$ and $t_{y}$ the subtree of $t_{x}$
starting at
$y,$ then $t_{y} \xbe T_{y}$.

(5) $<x,t_{x}> \xee <y,t_{y}>$ implies $y \xbe ran(f) \xcv \{x':rg(x'
)<rg(x)\}$ for some $f \xbe \xbP_{x}.$

$ \xcz $
\\[3ex]

\efa

Claim \ref{Claim 6.7} shows basic representation.

\bc

$\hspace{0.01em}$

% (+++*** Orig. No.:  Claim 6.7 )

\label{Claim 6.7}

$ \xcA U \xbe \xdy. \xbm (U)= \xbm_{ \xdz }(U)$

\ec

\subparagraph{
Proof:
}

$\hspace{0.01em}$

% (+++*** Orig.:  Proof: )

\label{Section Proof:}

By Claim \ref{Claim 6.1}, it suffices to show that for all $U \xbe
\xdy $
$x \xbe \xbm_{ \xdz }(U)$ $ \xcr $ $x \xbe U$ $ \xcu $ $ \xcE f \xbe
\xbP_{x}.ran(f) \xcs U= \xCQ.$

Fix $U \xbe \xdy.$

`` $ \xcp $ '': $x \xbe \xbm_{ \xdz }(U)$ $ \xcp $ ex. $<x,t_{x}>$ minimal
in $ \xdz \xex U,$ thus $x \xbe U$ and there is no $<y,t_{y}> \xbe \xdz,$
$<y,t_{y}> \xeb <x,t_{x}>,$ $y \xbe U.$ Let $f$ define the first part of
the set of children of
the root $x$ in $t_{x}$. If $ran(f) \xcs U \xEd \xCQ,$ if $y \xbe U$ is
a child of $x$ in $t_{x}$, and if
$t_{y}$ is the subtree of $t_{x}$ starting at $y,$ then $t_{y} \xbe T_{y}$
and $<y,t_{y}> \xeb <x,t_{x}>,$
contradicting minimality of $<x,t_{x}>$ in $ \xdz \xex U.$ So $ran(f) \xcs
U= \xCQ.$

`` $ \xcq $ '': Let $x \xbe U,$ and $ \xcE f \xbe \xbP_{x}.ran(f) \xcs U=
\xCQ.$ By Corollary \ref{Corollary 6.2},
there is no $x' \xbe U,$
$rg(x' )<rg(x).$ If $ \xdy_{x}= \xCQ,$ then the tree $tc_{x}$ has no $
\xem -$successors in $U,$ and
$<x,tc_{x}>$ is $ \xee -$minimal in $ \xdz \xex U.$ If $ \xdy_{x} \xEd
\xCQ $ and $f \xbe \xbP_{x}$ s.t. $ran(f) \xcs U= \xCQ,$ then
$<x,tf_{x}>$ is again $ \xee -$minimal in $ \xdz \xex U.$

$ \xcz $
\\[3ex]

We consider now the transitive closure of $ \xdz.$ (Recall that $
\xeb^{*}$ denotes the
transitive closure of $ \xeb.)$ Claim \ref{Claim 6.8} shows that
transitivity does not
destroy what we have achieved. The trees $tf_{x}$ play a crucial role in
the
demonstration.

\bc

$\hspace{0.01em}$

% (+++*** Orig. No.:  Claim 6.8 )

\label{Claim 6.8}

Let $ \xdz ' $ $:=$ $<$ $\{<x,t_{x}>:$ $x \xbe Z,$ $t_{x} \xbe T_{x}\}$,
$<x,t_{x}> \xee <y,t_{y}>$ iff $t_{x} \xem^{*}t_{y}$ $>.$ Then $ \xbm_{
\xdz }= \xbm_{ \xdz ' }.$

\ec

\subparagraph{
Proof:
}

$\hspace{0.01em}$

% (+++*** Orig.:  Proof: )

\label{Section Proof:}

Suppose there is $U \xbe \xdy,$ $x \xbe U,$ $x \xbe \xbm_{ \xdz }(U),$ $x
\xce \xbm_{ \xdz ' }(U).$
Then there must be an element $<x,t_{x}> \xbe \xdz $ with no $<x,t_{x}>
\xee <y,t_{y}>$ for any $y \xbe U.$
Let $f \xbe \xbP_{x}$ determine the first part of the set of children of
$x$ in $t_{x}$,
then $ran(f) \xcs U= \xCQ,$ consider $tf_{x}.$
All elements $w \xEd x$ of $tf_{x}$ are already in $ran(f),$ or
$rg(w)<rg(x)$ holds. (Note
that the elements chosen by rang in $tf_{x}$ continue by themselves or by
another
element of even smaller rang, but the rang order is transitive.) But all
$w$ s.t.
$rg(w)<rg(x)$ were already successors at level 1 of $x$ in $tf_{x}.$
By Corollary \ref{Corollary 6.2},
there is no $w \xbe U,$ $rg(w)<rg(x).$ Thus, no element $ \xEd x$ of
$tf_{x}$ is in $U.$
Thus there is no $<z,t_{z}> \xeb^{*}<x,tf_{x}>$ in $ \xdz $ with $z \xbe
U,$ so $<x,tf_{x}>$ is minimal
in $ \xdz ' \xex U,$ contradiction.

$ \xcz $ (Claim \ref{Claim 6.8} and Proposition \ref{Proposition 6.5})
\\[3ex]

% karl-search= End Hier-ARepr-General-Trans
\vspace{7mm}

% *************************************

\vspace{7mm}

% karl-search= End Hier-ARepr-General
\vspace{7mm}

% *************************************

\vspace{7mm}

%  6.4  &a-ranked smooth structures
%  6.4  &a-ranked smooth structures
% %
% =================================
\subsection{
$\xda$-ranked smooth structures
}

% {\LARGE karl-search= Start Hier-ARepr-Smooth }

\label{Section Hier-ARepr-Smooth}
\index{Section Hier-ARepr-Smooth}

% {\LARGE karl-search= Start Hier-ARepr-Smooth-Intro }

\label{Section Hier-ARepr-Smooth-Intro}
\index{Section Hier-ARepr-Smooth-Intro}

All smooth cases have a simple solution. We use one of our existing proofs
for
the not necessarily $ \xda -$ranked case, and add one litte result:

\bfa

$\hspace{0.01em}$

% (+++*** Orig. No.:  Fact Smooth-to-A )

\label{Fact Smooth-to-A}

Let $( \xbm \xda )$ hold, and let $ \xdz =< \xdx, \xeb >$ be a smooth
preferential structure
representing $ \xbm,$ i.e. $ \xbm = \xbm_{ \xdz }.$

Suppose that

$<x,i> \xeb <y,j>$ implies $rg(x) \xck rg(y).$

Define $ \xdz ':=< \xdx, \xer >$ where $<x,i> \xer <y,j>$ iff $<x,i>
\xeb <y,j>$ or $rg(x)<rg(y).$

Then $ \xdz ' $ is $ \xda -$ranked.

$ \xdz ' $ is smooth, too, and $ \xbm_{ \xdz }= \xbm_{ \xdz ' }=: \xbm '
.$

In addition, if $ \xeb $ is free from cycles, so is $ \xer,$ if $ \xeb $
is transitive, so is $ \xer.$

\efa

\subparagraph{
Proof:
}

$\hspace{0.01em}$

% (+++*** Orig.:  Proof: )

\label{Section Proof:}

$ \xda -$rankedness is trivial.

Suppose $<x,i>$ is $ \xeb -$minimal, but not $ \xer -$minimal. Then there
must be
$<y,j> \xer <x,i>,$ $<y,j> \xeB <x,i>,$ $y \xbe X,$ so
$rg(y)<rg(x).$ By $( \xbm \xda ),$ all $x \xbe \xbm (X)$ have minimal $
\xda -$rang among the elements of
$X,$ so this is impossible. Thus, $ \xbm -$minimal elements stay $ \xbm '
-$minimal, so
smoothness will also be preserved - remember that we increased the
relation.

By prerequisite, there cannot be any cycle involving only $ \xeb,$ but
the rang
order is free from cycles, too, and $ \xeb $ respects the rang order, so $
\xer $ is free
from cycles.

Let $ \xeb $ be transitive, so is the rang order. But if $<x,i> \xeb
<y,j>$ and
$rg(y)<rg(z)$ for some $<z,k>,$ then by prerequisite $rg(x) \xck rg(y),$
so $rg(x)<rg(z),$
so $<x,i> \xer <z,k>$ by definition. Likewise for $rg(x)<rg(y)$ and $<y,j>
\xeb <z,k>.$

$ \xcz $
\\[3ex]

All that remains to show then is that our constructions of smooth and of
smooth and transitive structures satisfy the condition

$<x,i> \xeb <y,j>$ implies $rg(x) \xck rg(y).$

% karl-search= End Hier-ARepr-Smooth-Intro
\vspace{7mm}

% *************************************

\vspace{7mm}

%  6.4.1  The basic smooth, not necessarily transitive case
%  6.4.1  The basic smooth, not necessarily transitive case
% %
% =========================================================
\subsubsection{
The basic smooth, not necessarily transitive case
}

% {\LARGE karl-search= Start Hier-ARepr-Smooth-Intrans }

\label{Section Hier-ARepr-Smooth-Intrans}
\index{Section Hier-ARepr-Smooth-Intrans}

% {\LARGE karl-search= Start Proposition A-Smooth-Complete }

\index{Proposition A-Smooth-Complete}

We will show here the following representation result:

\bp

$\hspace{0.01em}$

% (+++*** Orig. No.:  Proposition A-Smooth-Complete )

\label{Proposition A-Smooth-Complete}

Let $ \xda $ be given.

Let $ \xdy $ be closed under finite unions, and $ \xbm: \xdy \xcp \xdp
(Z).$
Then there is a $ \xdy -$smooth $ \xda -$ranked preferential structure $
\xdz,$ s.t. for all
$X \xbe \xdy $ $ \xbm (X)= \xbm_{ \xdz }(X)$ iff $ \xbm $ satisfies $(
\xbm \xcc ),$ $( \xbm PR),$ $( \xbm CUM),$ $( \xbm \xda ).$

% karl-search= End Proposition A-Smooth-Complete
\vspace{7mm}

% *************************************

\vspace{7mm}

\ep

To prove Proposition \ref{Proposition A-Smooth-Complete}, we first show and
prove:
\index{Proposition Smooth-Complete}

\bp

$\hspace{0.01em}$

% (+++*** Orig. No.:  Proposition Smooth-Complete )

\label{Proposition Smooth-Complete}

Let $ \xbm: \xdy \xcp \xdp (U)$ satisfy $( \xbm \xcc ),$ $( \xbm PR),$
and $( \xbm CUM),$ and the domain $ \xdy $ $( \xcv ).$

Then there is a $ \xdy -$smooth preferential structure $ \xdx $ s.t. $
\xbm = \xbm_{ \xdx }.$
\index{Proposition Smooth-Complete Proof}

\ep

\subparagraph{
Proof
}

$\hspace{0.01em}$

% (+++*** Orig.:  Proof )

% {\LARGE karl-search= Start Comment Smooth-Complete Proof }

\index{Comment Smooth-Complete Proof}

Outline: We first define a structure $ \xdz $ (in a way very similar to
Construction \ref{Construction Pref-Base})
which represents $ \xbm,$ but is not
necessarily $ \xdy -$smooth, refine
it to $ \xdz ' $ and show that $ \xdz ' $ represents $ \xbm $ too, and
that $ \xdz ' $ is $ \xdy -$smooth.

In the structure $ \xdz ',$ all pairs destroying smoothness in $ \xdz $
are successively
repaired, by adding minimal elements: If $<y,j>$ is not minimal, and has
no minimal
$<x,i>$ below it, we just add one such $<x,i>.$ As the repair process
might itself
generate such ``bad'' pairs, the process may have to be repeated infinitely
often.
Of course, one has to take care that the representation property is
preserved.

The proof given is close to the minimum one has to show (except that we
avoid
$H(U),$ instead of $U$ - as was done in the old proof of [Sch96-1]). We
could simplify
further, we do not, in order to stay closer to the construction that is
really
needed. The reader will find the simplification as building block of the
proof
of the transitive case. (In the simplified proof, we would consider for
$x,U$ s.t.
$x \xbe \xbm (U)$
the pairs $<x,g_{U}>$ with $g_{U} \xbe \xbP \{ \xbm (U \xcv Y):x \xbe Y
\xcC H(U)\},$ giving minimal elements.
For the $U$ s.t. $x \xbe U- \xbm (U),$ we would choose $<x,g>$ s.t. $g
\xbe \xbP \{ \xbm (Y):x \xbe Y \xbe \xdy \}$ with
$<x',g'_{U}> \xeb <x,g>$ for $<x',g'_{U}>$ as above.)

Construction \ref{Construction Smooth-Base} represents $ \xbm.$ The
structure will not yet
be smooth, we will mend it afterwards in Construction \ref{Construction
Smooth-Admiss}.

% karl-search= End Comment Smooth-Complete Proof
\vspace{7mm}

% *************************************

\vspace{7mm}

% {\LARGE karl-search= Start Comment HU }

\index{Comment HU}

\paragraph{
The constructions
}

$\hspace{0.01em}$

% (+++*** Orig.:  The constructions )

\label{Section The constructions}

$ \xdy $ will be closed under finite unions
throughout this Section. We first define $H(U),$ and show some facts
about it. $H(U)$ has an important role, for the
following
reason: If $u \xbe \xbm (U),$ but $u \xbe X- \xbm (X),$ then there is $x
\xbe \xbm (X)-H(U).$ Consequently,
to kill minimality of $u$ in $X,$ we can choose $x \xbe \xbm (X)-H(U),$ $x
\xeb u,$ without
interfering with u's minimality in $U.$ Moreover, if $x \xbe Y- \xbm (Y),$
then, by $x \xce H(U),$
$ \xbm (Y) \xcC H(U),$ so we can kill minimality of $x$ in $Y$ by choosing
some $y \xce H(U).$
Thus, even in the transitive case, we can leave $U$ to destroy minimality
of $u$ in
some $X,$ without ever having to come back into $U,$ it suffices to choose
sufficiently far from $U,$ i.e. outside $H(U).$ $H(U)$ is the right notion
of
``neighborhood''.

Note: Not all $z \xbe Z$ have to occur in our structure, therefore it is
quite
possible that $X \xbe \xdy,$ $X \xEd \xCQ,$ but $ \xbm_{ \xdz }(X)= \xCQ
.$ This is why we have introduced
the set $K$ in Definition \ref{Definition K} and such $X$ will be subsets
of $Z-$K.

Let now $ \xbm: \xdy \xcp \xdp (Z).$

% karl-search= End Comment HU
\vspace{7mm}

% *************************************

\vspace{7mm}

% {\LARGE karl-search= Start Definition HU }

\index{Definition HU}

\bd

$\hspace{0.01em}$

% (+++*** Orig. No.:  Definition HU )

\label{Definition HU}

Define $H(U)$ $:=$ $ \xcV \{X: \xbm (X) \xcc U\}.$

% karl-search= End Definition HU
\vspace{7mm}

% *************************************

\vspace{7mm}

% {\LARGE karl-search= Start Definition K }

\index{Definition K}

\ed

\bd

$\hspace{0.01em}$

% (+++*** Orig. No.:  Definition K )

\label{Definition K}

Let $K:=\{x \xbe Z:$ $ \xcE X \xbe \xdy.x \xbe \xbm (X)\}$

% karl-search= End Definition K
\vspace{7mm}

% *************************************

\vspace{7mm}

% {\LARGE karl-search= Start Fact HU-1 }

\index{Fact HU-1}

\ed

\bfa

$\hspace{0.01em}$

% (+++*** Orig. No.:  Fact HU-1 )

\label{Fact HU-1}

$( \xbm \xcc )$ $+$ $( \xbm PR)$ $+$ $( \xbm CUM)$ $+$ $( \xcv )$ entail:

(1) $ \xbm (A) \xcc B$ $ \xcp $ $ \xbm (A \xcv B)= \xbm (B)$

(2) $ \xbm (X) \xcc U,$ $U \xcc Y$ $ \xcp $ $ \xbm (Y \xcv X)= \xbm (Y)$

(3) $ \xbm (X) \xcc U,$ $U \xcc Y$ $ \xcp $ $ \xbm (Y) \xcs X \xcc \xbm
(U)$

(4) $ \xbm (X) \xcc U$ $ \xcp $ $ \xbm (U) \xcs X \xcc \xbm (X)$

(5) $U \xcc A,$ $ \xbm (A) \xcc H(U)$ $ \xcp $ $ \xbm (A) \xcc U$

(6) Let $x \xbe K,$ $Y \xbe \xdy,$ $x \xbe Y- \xbm (Y),$ then $ \xbm (Y)
\xEd \xCQ.$

% karl-search= End Fact HU-1
\vspace{7mm}

% *************************************

\vspace{7mm}

% {\LARGE karl-search= Start Fact HU-1 Proof }

\index{Fact HU-1 Proof}

\efa

\subparagraph{
Proof:
}

$\hspace{0.01em}$

% (+++*** Orig.:  Proof: )

\label{Section Proof:}

(1) $ \xbm (A) \xcc B$ $ \xcp $ $ \xbm (A \xcv B) \xcc \xbm (A) \xcv \xbm
(B) \xcc B$ $ \xcp_{( \xbm CUM)}$ $ \xbm (B)= \xbm (A \xcv B).$

(2) trivial by (1).

(3) $ \xbm (Y) \xcs X$ $=$ (by (2)) $ \xbm (Y \xcv X) \xcs X$ $ \xcc $ $
\xbm (Y \xcv X) \xcs (X \xcv U)$ $ \xcc $ (by $( \xbm PR))$
$ \xbm (X \xcv U)$ $=$ (by (1)) $ \xbm (U).$

(4) $ \xbm (U) \xcs X$ $=$ $ \xbm (X \xcv U) \xcs X$ by (1) $ \xcc $ $
\xbm (X)$ by $( \xbm PR)$

(5) Let $U \xcc A,$ $ \xbm (A) \xcc H(U).$ So $ \xbm (A)$ $=$ $ \xcV \{
\xbm (A) \xcs Y: \xbm (Y) \xcc U\}$ $ \xcc $
$ \xbm (U)$ $ \xcc $ $U$ by (3).

(6) Suppose $x \xbe \xbm (X),$ $ \xbm (Y)= \xCQ $ $ \xcp $ $ \xbm (Y) \xcc
X,$ so by (4) $Y \xcs \xbm (X) \xcc \xbm (Y),$ so
$x \xbe \xbm (Y).$

$ \xcz $
\\[3ex]

% karl-search= End Fact HU-1 Proof
\vspace{7mm}

% *************************************

\vspace{7mm}

% {\LARGE karl-search= Start Fact HU-2 }

\index{Fact HU-2}

The following Fact \ref{Fact HU-2} contains the basic properties
of $ \xbm $ and $H(U)$ which we will need for the representation
construction.

\bfa

$\hspace{0.01em}$

% (+++*** Orig. No.:  Fact HU-2 )

\label{Fact HU-2}

Let A, $U,$ $U',$ $Y$ and all $A_{i}$ be in $ \xdy.$ Let $( \xbm \xcc )$
$+$ $( \xbm PR)$ $+$ $( \xcv )$ hold.

(1) $A= \xcV \{A_{i}:i \xbe I\}$ $ \xcp $ $ \xbm (A) \xcc \xcV \{ \xbm
(A_{i}):i \xbe I\},$

(2) $U \xcc H(U),$ and $U \xcc U' \xcp H(U) \xcc H(U' ),$

(3) $ \xbm (U \xcv Y)-H(U) \xcc \xbm (Y).$

If, in addition, $( \xbm CUM)$ holds, then we also have:

(4) $U \xcc A,$ $ \xbm (A) \xcc H(U)$ $ \xcp $ $ \xbm (A) \xcc U,$

(5) $ \xbm (Y) \xcc H(U)$ $ \xcp $ $Y \xcc H(U)$ and $ \xbm (U \xcv Y)=
\xbm (U),$

(6) $x \xbe \xbm (U),$ $x \xbe Y- \xbm (Y)$ $ \xcp $ $Y \xcC H(U),$

(7) $Y \xcC H(U)$ $ \xcp $ $ \xbm (U \xcv Y) \xcC H(U).$

% karl-search= End Fact HU-2
\vspace{7mm}

% *************************************

\vspace{7mm}

% {\LARGE karl-search= Start Fact HU-2 Proof }

\index{Fact HU-2 Proof}

\efa

\subparagraph{
Proof:
}

$\hspace{0.01em}$

% (+++*** Orig.:  Proof: )

\label{Section Proof:}

(1) $ \xbm (A) \xcs A_{j} \xcc \xbm (A_{j}) \xcc \xcV \xbm (A_{i}),$ so by
$ \xbm (A) \xcc A= \xcV A_{i}$ $ \xbm (A) \xcc \xcV \xbm (A_{i}).$

(2) trivial.

(3) $ \xbm (U \xcv Y)-H(U)$ $ \xcc_{(2)}$ $ \xbm (U \xcv Y)-U$ $ \xcc_{(
\xbm \xcc )}$ $ \xbm (U \xcv Y) \xcs Y$ $ \xcc_{( \xbm PR)}$ $ \xbm (Y).$

(4) This is Fact \ref{Fact HU-1} (5).

(5) Let $ \xbm (Y) \xcc H(U),$ then by $ \xbm (U) \xcc H(U)$ and (1) $
\xbm (U \xcv Y) \xcc \xbm (U) \xcv \xbm (Y) \xcc H(U),$
so by (4) $ \xbm (U \xcv Y) \xcc U$ and $U \xcv Y \xcc H(U).$ Moreover, $
\xbm (U \xcv Y) \xcc U \xcc U \xcv Y$ $ \xcp_{( \xbm CUM)}$ $ \xbm (U \xcv
Y)= \xbm (U).$

(6) If not, $Y \xcc H(U),$ so $ \xbm (Y) \xcc H(U),$ so $ \xbm (U \xcv Y)=
\xbm (U)$ by (5), but $x \xbe Y- \xbm (Y)$ $ \xcp_{( \xbm PR)}$
$x \xce \xbm (U \xcv Y)= \xbm (U),$ $contradiction.$

(7) $ \xbm (U \xcv Y) \xcc H(U)$ $ \xcp_{(5)}$ $U \xcv Y \xcc H(U).$
$ \xcz $
\\[3ex]

% karl-search= End Fact HU-2 Proof
\vspace{7mm}

% *************************************

\vspace{7mm}

% {\LARGE karl-search= Start Definition Gamma-x }

\index{Definition Gamma-x}

\bd

$\hspace{0.01em}$

% (+++*** Orig. No.:  Definition Gamma-x )

\label{Definition Gamma-x}

For $x \xbe Z,$ let $ \xdw_{x}:=\{ \xbm (Y)$: $Y \xbe \xdy $ $ \xcu $ $x
\xbe Y- \xbm (Y)\},$ $ \xbG_{x}:= \xbP \xdw_{x}$, and $K:=\{x \xbe Z$: $
\xcE X \xbe \xdy.x \xbe \xbm (X)\}.$

\ed

Note that we consider here now $ \xbm (Y)$ in $ \xdw_{x}$, and not $Y$ as
in $ \xdy_{x}$
in Definition \ref{Definition Y-Pi-x}.

% karl-search= End Definition Gamma-x
\vspace{7mm}

% *************************************

\vspace{7mm}

% {\LARGE karl-search= Start Remark Gamma-x }

\index{Remark Gamma-x}

\br

$\hspace{0.01em}$

% (+++*** Orig. No.:  Remark Gamma-x )

\label{Remark Gamma-x}

Assume now $( \xbm \xcc ),$ $( \xbm PR),$ $( \xbm CUM),$ $( \xcv )$ to
hold.

(1) $x \xbe K$ $ \xcp $ $ \xbG_{x} \xEd \xCQ,$

(2) $g \xbe \xbG_{x}$ $ \xcp $ $ran(g) \xcc K.$

% karl-search= End Remark Gamma-x
\vspace{7mm}

% *************************************

\vspace{7mm}

% {\LARGE karl-search= Start Remark Gamma-x Proof }

\index{Remark Gamma-x Proof}

\er

\subparagraph{
Proof:
}

$\hspace{0.01em}$

% (+++*** Orig.:  Proof: )

\label{Section Proof:}

(1) We have to show that $Y \xbe \xdy,$ $x \xbe Y- \xbm (Y)$ $ \xcp $ $
\xbm (Y) \xEd \xCQ.$ This was shown
in Fact \ref{Fact HU-1} (6).

(2) By definition, $ \xbm (Y) \xcc K$ for all $Y \xbe \xdy.$
$ \xcz $
\\[3ex]

% karl-search= End Remark Gamma-x Proof
\vspace{7mm}

% *************************************

\vspace{7mm}

% {\LARGE karl-search= Start Claim Cum-Mu-f }

\index{Claim Cum-Mu-f}

The following claim is the analogue of Claim \ref{Claim Mu-f} above.

\bc

$\hspace{0.01em}$

% (+++*** Orig. No.:  Claim Cum-Mu-f )

\label{Claim Cum-Mu-f}

Assume now $( \xbm \xcc ),$ $( \xbm PR),$ $( \xbm CUM),$ $( \xcv )$ to
hold.

Let $U \xbe \xdy,$ $x \xbe K.$ Then

(1) $x \xbe \xbm (U)$ $ \xcr $ $x \xbe U$ $ \xcu $ $ \xcE f \xbe
\xbG_{x}.ran(f) \xcs U= \xCQ,$

(2) $x \xbe \xbm (U)$ $ \xcr $ $x \xbe U$ $ \xcu $ $ \xcE f \xbe
\xbG_{x}.ran(f) \xcs H(U)= \xCQ.$

% karl-search= End Claim Cum-Mu-f
\vspace{7mm}

% *************************************

\vspace{7mm}

% {\LARGE karl-search= Start Claim Cum-Mu-f Proof }

\index{Claim Cum-Mu-f Proof}

\ec

\subparagraph{
Proof:
}

$\hspace{0.01em}$

% (+++*** Orig.:  Proof: )

\label{Section Proof:}

(1)
Case 1: $ \xdw_{x}= \xCQ,$ thus $ \xbG_{x}=\{ \xCQ \}.$
`` $ \xcp $ '': Take $f:= \xCQ.$
`` $ \xcq $ '': $x \xbe U \xbe \xdy,$ $ \xdw_{x}= \xCQ $ $ \xcp $ $x \xbe
\xbm (U)$ by definition of $ \xdw_{x}.$

Case 2: $ \xdw_{x} \xEd \xCQ.$
`` $ \xcp $ '': Let $x \xbe \xbm (U) \xcc U.$ It suffices to show $Y \xbe
\xdw_{x}$ $ \xcp $ $ \xbm (Y)-H(U) \xEd \xCQ.$
But $Y \xbe \xdw_{x}$ $ \xcp $ $x \xbe Y- \xbm (Y)$ $ \xcp $ (by Fact \ref{Fact
HU-2}, (6))
$Y \xcC H(U)$ $ \xcp $ (by Fact \ref{Fact HU-2}, (5)) $ \xbm (Y) \xcC
H(U).$
`` $ \xcq $ '': If $x \xbe U- \xbm (U),$ $U \xbe \xdw_{x}$, moreover $
\xbG_{x} \xEd \xCQ $ by
Remark \ref{Remark Gamma-x}, (1) and thus (or
by the same argument) $ \xbm (U) \xEd \xCQ,$ so $ \xcA f \xbe
\xbG_{x}.ran(f) \xcs U \xEd \xCQ.$

(2): The proof is verbatim the same as for (1).
$ \xcz $
\\[3ex]

% karl-search= End Claim Cum-Mu-f Proof
\vspace{7mm}

% *************************************

\vspace{7mm}

\subparagraph{
Proof: (Prop. 6.14)
}

$\hspace{0.01em}$

% (+++*** Orig.:  Proof: (Prop. 6.14) )

\label{Section Proof: (Prop. 6.14)}

% {\LARGE karl-search= Start Construction Smooth-Base }

\index{Construction Smooth-Base}

\bcs

$\hspace{0.01em}$

% (+++*** Orig. No.:  Construction Smooth-Base )

\label{Construction Smooth-Base}

(Construction of $ \xdz )$
Let $ \xdx $ $:=$ $\{<x,g>$: $x \xbe K,$ $g \xbe \xbG_{x}\},$ $<x',g' >
\xeb <x,g>$ $: \xcr $ $x' \xbe ran(g),$

$ \xdz:=< \xdx, \xeb >.$

% karl-search= End Construction Smooth-Base
\vspace{7mm}

% *************************************

\vspace{7mm}

% {\LARGE karl-search= Start Claim Smooth-Base }

\index{Claim Smooth-Base}

\ecs

\bc

$\hspace{0.01em}$

% (+++*** Orig. No.:  Claim Smooth-Base )

\label{Claim Smooth-Base}

$ \xcA U \xbe \xdy. \xbm (U)= \xbm_{ \xdz }(U)$

% karl-search= End Claim Smooth-Base
\vspace{7mm}

% *************************************

\vspace{7mm}

% {\LARGE karl-search= Start Claim Smooth-Base Proof }

\index{Claim Smooth-Base Proof}

\ec

\subparagraph{
Proof:
}

$\hspace{0.01em}$

% (+++*** Orig.:  Proof: )

\label{Section Proof:}

Case 1: $x \xce K.$ Then $x \xce \xbm (U)$ and $x \xce \xbm_{ \xdz }(U).$

Case 2: $x \xbe K.$
By Claim \ref{Claim Cum-Mu-f}, (1) it suffices to show that for all $U
\xbe \xdy $
$x \xbe \xbm_{ \xdz }(U)$ $ \xcr $ $x \xbe U$ $ \xcu $ $ \xcE f \xbe
\xbG_{x}.ran(f) \xcs U= \xCQ.$
Fix $U \xbe \xdy.$
`` $ \xcp $ '': $x \xbe \xbm_{ \xdz }(U)$ $ \xcp $ ex. $<x,f>$ minimal in
$ \xdx \xex U,$ thus $x \xbe U$ and there is no
$<x',f' > \xeb <x,f>,$ $x' \xbe U,$ $x' \xbe K.$ But if $x' \xbe K,$ then
by
Remark \ref{Remark Gamma-x}, (1), $ \xbG_{x' } \xEd \xCQ,$
so we find suitable $f'.$ Thus, $ \xcA x' \xbe ran(f).x' \xce U$ or $x'
\xce K.$ But $ran(f) \xcc K,$ so
$ran(f) \xcs U= \xCQ.$
`` $ \xcq $ '': If $x \xbe U,$ $f \xbe \xbG_{x}$ s.t. $ran(f) \xcs U= \xCQ
,$ then $<x,f>$ is minimal in $ \xdx \xex U.$

$ \xcz $
\\[3ex]

% karl-search= End Claim Smooth-Base Proof
\vspace{7mm}

% *************************************

\vspace{7mm}

% {\LARGE karl-search= Start Construction Smooth-Admiss }

\index{Construction Smooth-Admiss}

We now construct the refined structure $ \xdz '.$

\bcs

$\hspace{0.01em}$

% (+++*** Orig. No.:  Construction Smooth-Admiss )

\label{Construction Smooth-Admiss}

(Construction of $ \xdz ' )$

$ \xbs $ is called $x-$admissible sequence iff

1. $ \xbs $ is a sequence of length $ \xck \xbo,$ $ \xbs =\{ \xbs_{i}:i
\xbe \xbo \},$

2. $ \xbs_{o} \xbe \xbP \{ \xbm (Y)$: $Y \xbe \xdy $ $ \xcu $ $x \xbe Y-
\xbm (Y)\},$

3. $ \xbs_{i+1} \xbe \xbP \{ \xbm (X)$: $X \xbe \xdy $ $ \xcu $ $x \xbe
\xbm (X)$ $ \xcu $ $ran( \xbs_{i}) \xcs X \xEd \xCQ \}.$

By 2., $ \xbs_{0}$ minimizes $x,$ and by 3., if $x \xbe \xbm (X),$ and
$ran( \xbs_{i}) \xcs X \xEd \xCQ,$ i.e. we
have destroyed minimality of $x$ in $X,$ $x$ will be above some $y$
minimal in $X$ to
preserve smoothness.

Let $ \xbS_{x}$ be the set of $x-$admissible sequences, for $ \xbs \xbe
\xbS_{x}$ let $ \wt{ \xbs }:= \xcV \{ran( \xbs_{i}):i \xbe \xbo \}.$
Note that by the argument in the proof of
Remark \ref{Remark Gamma-x}, (1),
$ \xbS_{x} \xEd \xCQ,$ if $x \xbe K.$

Let $ \xdx ' $ $:=$ $\{<x, \xbs >$: $x \xbe K$ $ \xcu $ $ \xbs \xbe
\xbS_{x}\}$ and $<x', \xbs ' > \xeb ' <x, \xbs >$ $: \xcr $ $x' \xbe \wt{
\xbs }$.
Finally, let $ \xdz ':=< \xdx ', \xeb ' >,$ and $ \xbm ':= \xbm_{ \xdz
' }.$

\ecs

It is now easy to show that $ \xdz ' $ represents $ \xbm,$ and that $
\xdz ' $ is smooth.
For $x \xbe \xbm (U),$ we construct a special $x-$admissible sequence $
\xbs^{x,U}$ using the
properties of $H(U).$

% karl-search= End Construction Smooth-Admiss
\vspace{7mm}

% *************************************

\vspace{7mm}

% {\LARGE karl-search= Start Claim Smooth-Admiss-1 }

\index{Claim Smooth-Admiss-1}

\bc

$\hspace{0.01em}$

% (+++*** Orig. No.:  Claim Smooth-Admiss-1 )

\label{Claim Smooth-Admiss-1}

For all $U \xbe \xdy $ $ \xbm (U)= \xbm_{ \xdz }(U)= \xbm ' (U).$

% karl-search= End Claim Smooth-Admiss-1
\vspace{7mm}

% *************************************

\vspace{7mm}

% {\LARGE karl-search= Start Claim Smooth-Admiss-1 Proof }

\index{Claim Smooth-Admiss-1 Proof}

\ec

\subparagraph{
Proof:
}

$\hspace{0.01em}$

% (+++*** Orig.:  Proof: )

\label{Section Proof:}

If $x \xce K,$ then $x \xce \xbm_{ \xdz }(U),$ and $x \xce \xbm ' (U)$ for
any $U.$ So assume $x \xbe K.$ If $x \xbe U$ and
$x \xce \xbm_{ \xdz }(U),$ then for all $<x,f> \xbe \xdx,$ there is $<x'
,f' > \xbe \xdx $ with $<x',f' > \xeb <x,f>$ and
$x' \xbe U.$ Let now $<x, \xbs > \xbe \xdx ',$ then $<x, \xbs_{0}> \xbe
\xdx,$ and let $<x',f' > \xeb <x, \xbs_{0}>$ in $ \xdz $ with
$x' \xbe U.$ As $x' \xbe K,$ $ \xbS_{x' } \xEd \xCQ,$ let $ \xbs ' \xbe
\xbS_{x' }$. Then $<x', \xbs ' > \xeb ' <x, \xbs >$ in $ \xdz '.$ Thus
$x \xce \xbm ' (U).$
Thus, for all $U \xbe \xdy,$ $ \xbm ' (U) \xcc \xbm_{ \xdz }(U)= \xbm
(U).$

It remains to show $x \xbe \xbm (U) \xcp x \xbe \xbm ' (U).$

Assume $x \xbe \xbm (U)$ (so $x \xbe K),$ $U \xbe \xdy,$ we will
construct minimal $ \xbs,$ i.e. show that
there is $ \xbs^{x,U} \xbe \xbS_{x}$ s.t. $ \wt{ \xbs^{x,U}} \xcs U= \xCQ
.$ We construct this $ \xbs^{x,U}$ inductively, with the
stronger property that $ran( \xbs^{x,U}_{i}) \xcs H(U)= \xCQ $ for all $i
\xbe \xbo.$

$ \ul{ \xbs^{x,U}_{0}:}$
$x \xbe \xbm (U),$ $x \xbe Y- \xbm (Y)$ $ \xcp $ $ \xbm (Y)-H(U) \xEd \xCQ
$ by Fact \ref{Fact HU-2}, $(6)+(5).$
Let $ \xbs^{x,U}_{0}$ $ \xbe $ $ \xbP \{ \xbm (Y)-H(U):$ $Y \xbe \xdy,$
$x \xbe Y- \xbm (Y)\},$ so $ran( \xbs^{x,U}_{0}) \xcs H(U)= \xCQ.$

$ \ul{ \xbs^{x,U}_{i} \xcp \xbs^{x,U}_{i+1}:}$
By induction hypothesis, $ran( \xbs^{x,U}_{i}) \xcs H(U)= \xCQ.$ Let $X
\xbe \xdy $ be s.t. $x \xbe \xbm (X),$
$ran( \xbs^{x,U}_{i}) \xcs X \xEd \xCQ.$ Thus $X \xcC H(U),$ so $ \xbm (U
\xcv X)-H(U) \xEd \xCQ $ by Fact \ref{Fact HU-2}, (7).
Let $ \xbs^{x,U}_{i+1}$ $ \xbe $ $ \xbP \{ \xbm (U \xcv X)-H(U):$ $X \xbe
\xdy,$ $x \xbe \xbm (X),$ $ran( \xbs^{x,U}_{i}) \xcs X \xEd \xCQ \},$ so
$ran( \xbs^{x,U}_{i+1}) \xcs H(U)= \xCQ.$
As $ \xbm (U \xcv X)-H(U) \xcc \xbm (X)$ by Fact \ref{Fact HU-2},
(3), the construction satisfies the
$x-$admissibility condition.
$ \xcz $
\\[3ex]

% karl-search= End Claim Smooth-Admiss-1 Proof
\vspace{7mm}

% *************************************

\vspace{7mm}

% {\LARGE karl-search= Start Claim Smooth-Admiss-2 }

\index{Claim Smooth-Admiss-2}

We now show:

\bc

$\hspace{0.01em}$

% (+++*** Orig. No.:  Claim Smooth-Admiss-2 )

\label{Claim Smooth-Admiss-2}

$ \xdz ' $ is $ \xdy -$smooth.

% karl-search= End Claim Smooth-Admiss-2
\vspace{7mm}

% *************************************

\vspace{7mm}

% {\LARGE karl-search= Start Claim Smooth-Admiss-2 Proof }

\index{Claim Smooth-Admiss-2 Proof}

\ec

\subparagraph{
Proof:
}

$\hspace{0.01em}$

% (+++*** Orig.:  Proof: )

\label{Section Proof:}

Let $X \xbe \xdy,$ $<x, \xbs > \xbe \xdx ' \xex X.$

Case 1, $x \xbe X- \xbm (X):$ Then $ran( \xbs_{0}) \xcs \xbm (X) \xEd \xCQ
,$ let $x' \xbe ran( \xbs_{0}) \xcs \xbm (X).$ Moreover,
$ \xbm (X) \xcc K.$ Then for all $<x', \xbs ' > \xbe \xdx ' $ $<x', \xbs
' > \xeb <x, \xbs >.$ But $<x', \xbs^{x',X}>$ as
constructed in the proof of Claim \ref{Claim Smooth-Admiss-1}
is minimal in $ \xdx ' \xex X.$

Case 2, $x \xbe \xbm (X)= \xbm_{ \xdz }(X)= \xbm ' (X):$ If $<x, \xbs >$
is minimal in $ \xdx ' \xex X,$ we are done.
So suppose there is $<x', \xbs ' > \xeb <x, \xbs >,$ $x' \xbe X.$ Thus
$x' \xbe \wt{ \xbs }.$ Let
$x' \xbe ran( \xbs_{i}).$ So $x \xbe \xbm (X)$ and $ran( \xbs_{i}) \xcs X
\xEd \xCQ.$ But
$ \xbs_{i+1} \xbe \xbP \{ \xbm (X' )$: $X' \xbe \xdy $ $ \xcu $ $x \xbe
\xbm (X' )$ $ \xcu $ $ran( \xbs_{i}) \xcs X' \xEd \xCQ \},$ so $X$ is one
of the $X',$
moreover $ \xbm (X) \xcc K,$ so there is $x'' \xbe \xbm (X) \xcs ran(
\xbs_{i+1}) \xcs K,$ so for all $<x'', \xbs '' > \xbe \xdx ' $
$<x'', \xbs '' > \xeb <x, \xbs >.$ But again $<x'', \xbs^{x'',X}>$ as
constructed in the proof of
Claim \ref{Claim Smooth-Admiss-1} is minimal in $ \xdx ' \xex X.$

$ \xcz $
\\[3ex]

% karl-search= End Claim Smooth-Admiss-2 Proof
\vspace{7mm}

% *************************************

\vspace{7mm}

% {\LARGE karl-search= Start Proposition A-Smooth-Complete Proof }

\index{Proposition A-Smooth-Complete Proof}

\paragraph{
Proof of Proposition \ref{Proposition A-Smooth-Complete}
}

$\hspace{0.01em}$

% (+++*** Orig.:  Proof of  {Proposition A-Smooth-Complete } )

\label{Section Proof of  {Proposition A-Smooth-Complete }}

Consider the construction in the proof of Proposition \ref{Proposition
Smooth-Complete}.
We have to show that it respects the rang order with respect
to $ \xda,$ i.e. that $<x', \xbs ' > \xeb ' <x, \xbs >$ implies $rg(x' )
\xck rg(x).$ This is easy: By
definition, $x' \xbe \xcV \{ran( \xbs_{i}):i \xbe \xbo \}.$
If $x' \xbe ran( \xbs_{0}),$ then for some $Y$ $x' \xbe \xbm (Y),$ $x \xbe
Y- \xbm (Y),$ so $rg(x' ) \xck rg(x)$
by $( \xbm \xda ).$
If $x' \xbe ran( \xbs_{i}),$ $i>0,$ then for some $X$ $x',x \xbe \xbm
(X),$ so $rg(x)=rg(x' )$ by $( \xbm \xda ).$

$ \xcz $ (Proposition \ref{Proposition A-Smooth-Complete})
\\[3ex]

% karl-search= End Proposition A-Smooth-Complete Proof
\vspace{7mm}

% *************************************

\vspace{7mm}

% karl-search= End Hier-ARepr-Smooth-Intrans
\vspace{7mm}

% *************************************

\vspace{7mm}

%  6.4.2  The transitive smooth case
%  6.4.2  The transitive smooth case
% %
% ==================================
\subsubsection{
The transitive smooth case
}

% {\LARGE karl-search= Start Hier-ARepr-Smooth-Trans }

\label{Section Hier-ARepr-Smooth-Trans}
\index{Section Hier-ARepr-Smooth-Trans}

% {\LARGE karl-search= Start Proposition A-Smooth-Trans-Complete }

\index{Proposition A-Smooth-Trans-Complete}

We will show here the transitive analogon of
Proposition \ref{Proposition A-Smooth-Complete}:

\bp

$\hspace{0.01em}$

% (+++*** Orig. No.:  Proposition A-Smooth-Trans-Complete )

\label{Proposition A-Smooth-Trans-Complete}

Let $ \xda $ be given.

Let $ \xdy $ be closed under finite unions, and $ \xbm: \xdy \xcp \xdp
(Z).$
Then there is a $ \xdy -$smooth $ \xda -$ranked transitive preferential
structure $ \xdz,$ s.t.
for all
$X \xbe \xdy $ $ \xbm (X)= \xbm_{ \xdz }(X)$ iff $ \xbm $ satisfies $(
\xbm \xcc ),$ $( \xbm PR),$ $( \xbm CUM),$ $( \xbm \xda ).$

% karl-search= End Proposition A-Smooth-Trans-Complete
\vspace{7mm}

% *************************************

\vspace{7mm}

\ep

To prove Proposition \ref{Proposition A-Smooth-Trans-Complete}, we
first show and prove:
\index{Proposition Smooth-Complete-Trans}

\bp

$\hspace{0.01em}$

% (+++*** Orig. No.:  Proposition Smooth-Complete-Trans )

\label{Proposition Smooth-Complete-Trans}

Let $ \xbm: \xdy \xcp \xdp (U)$ satisfy $( \xbm \xcc ),$ $( \xbm PR),$
and $( \xbm CUM),$ and the domain $ \xdy $ $( \xcv ).$

Then there is a transitive $ \xdy -$smooth preferential structure $ \xdx $
s.t. $ \xbm = \xbm_{ \xdx }.$
\index{Proposition Smooth-Complete-Trans Proof}

\ep

\subparagraph{
Proof
}

$\hspace{0.01em}$

% (+++*** Orig.:  Proof )

% {\LARGE karl-search= Start Discussion Smooth-Trans }

\index{Discussion Smooth-Trans}

\paragraph{
Discussion Smooth-Trans
}

$\hspace{0.01em}$

% (+++*** Orig.:  Discussion Smooth-Trans )

\label{Section Discussion Smooth-Trans}

In a certain way, it is not surprising that transitivity does not impose
stronger conditions in the smooth case either. Smoothness is itself a weak
kind of
transitivity: If an element is not minimal, then there is a minimal
element
below it, i.e., $x \xee y$ with $y$ not minimal is possible, there is $z'
\xeb y,$ but
then there is $z$ minimal with $x \xee z.$ This is ``almost'' $x \xee z',$
transitivity.

To obtain representation,
we will combine here the ideas of the smooth, but not necessarily
transitive
case with those of the general transitive case - as the reader will have
suspected. Thus, we will index again with trees, and work with (suitably
adapted) admissible sequences for the construction of the trees. In the
construction of the admissible sequences, we were careful to repair all
damage
done in previous steps. We have to add now reparation of all damage done
by
using transitivity, i.e., the transitivity of the relation might destroy
minimality, and we have to construct minimal elements below all elements
for
which we thus destroyed minimality. Both cases are combined by considering
immediately all $Y$ s.t. $x \xbe Y-H(U).$ Of course, the properties
described in
Fact \ref{Fact HU-2} play again a central role.

The (somewhat complicated) construction will be commented on in more
detail
below.

Note that even beyond Fact \ref{Fact HU-2}, closure of the domain
under finite
unions is used in the construction of the trees. This - or something like
it -
is necessary, as we have to respect the hulls of all elements treated so
far
(the predecessors), and not only of the first element, because of
transitivity.
For the same reason, we need more bookkeeping, to annotate all the hulls
(or
the union of the respective $U' $s) of all predecessors to be respected.

To summarize: we combine the ideas from the transitive general case and
the
simple smooth case, using the crucial Fact \ref{Fact HU-2}
to show that the construction
goes through. The construction leaves still some freedom, and
modifications
are possible as indicated below in the course of the proof.

Recall that $ \xdy $ will be closed under finite unions
in this section, and let again $ \xbm: \xdy \xcp \xdp (Z).$

We have to adapt Construction \ref{Construction Smooth-Admiss}
(x-admissible sequences) to the
transitive situation, and to our construction with trees. If $< \xCQ,x>$
is the root,
$ \xbs_{0} \xbe \xbP \{ \xbm (Y):x \xbe Y- \xbm (Y)\}$ determines some
children of the root.
To preserve smoothness, we have to compensate and add
other children by the $ \xbs_{i+1}:$ $ \xbs_{i+1} \xbe \xbP \{ \xbm (X):x
\xbe \xbm (X),$ $ran( \xbs_{i}) \xcs X \xEd \xCQ \}.$
On the other hand, we have to pursue the same construction for the
children so
constructed. Moreover, these indirect children have to be added to those
children
of the root, which have to be compensated (as the first children are
compensated
by $ \xbs_{1})$ to preserve smoothness. Thus, we build the tree in a
simultaneous vertical
and horizontal induction.

This construction can be simplified, by considering immediately all $Y
\xbe \xdy $ s.t.
$x \xbe Y \xcC H(U)$ - independent of whether $x \xce \xbm (Y)$ (as done
in $ \xbs_{0}),$ or whether
$x \xbe \xbm (Y),$ and some child $y$ constructed before is in $Y$ (as
done in the $ \xbs_{i+1}),$ or
whether $x \xbe \xbm (Y),$ and some indirect child $y$ of $x$ is in $Y$
(to take care of
transitivity, as indicated above). We make this simplified construction.

There are two ways to proceed. First, we can take as $ \xej^{*}$ in the
trees
the transitive closure of $ \xej.$ Second, we can deviate from the idea
that
children are chosen by selection functions $f,$ and take nonempty subsets
of
elements instead, making more elements children than in the first case. We
take
the first alternative, as it is more in the spirit of the construction.

We will suppose for simplicity that $Z=K$ - the general case in easy to
obtain,
but complicates the picture.

For each $x \xbe Z,$ we construct trees $t_{x}$, which will be used to
index
different copies of $x,$ and control the relation $ \xeb.$

These trees $t_{x}$ will have the following form:

(a) the root of $t$ is $< \xCQ,x>$ or $<U,x>$ with $U \xbe \xdy $ and $x
\xbe \xbm (U),$

(b) all other nodes are pairs $<Y,y>,$ $Y \xbe \xdy,$ $y \xbe \xbm (Y),$

(c) $ht(t) \xck \xbo,$

(d) if $<Y,y>$ is an element in $t_{x},$ then there is some $ \xdy (y)
\xcc \{W \xbe \xdy:y \xbe W\},$ and
$f \xbe \xbP \{ \xbm (W):W \xbe \xdy (y)\}$ s.t. the set of children of
$<Y,y>$ is $\{<Y \xcv W,f(W)>:$ $W \xbe \xdy (y)\}.$

The first coordinate is used for bookkeeping when constructing children,
in
particular for condition (d).

The relation $ \xeb $ will essentially be determined by the subtree
relation.

We first construct the trees $t_{x}$ for those sets $U$ where $x \xbe \xbm
(U),$ and then take
care of the others. In the construction for the minimal elements,
at each level $n>0,$ we may have several ways to choose a selection
function $f_{n}$,
and each such choice leads to the construction of a different tree - we
construct all these trees. (We could also construct only one tree, but
then
the choice would have to be made coherently for different $x,U.$ It is
simpler to
construct more trees than necessary.)

We control the relation by indexing with trees, just as it was done in the
not
necessarily smooth case before.

% karl-search= End Discussion Smooth-Trans
\vspace{7mm}

% *************************************

\vspace{7mm}

% {\LARGE karl-search= Start Definition Tree-TC }

\index{Definition Tree-TC}

\bd

$\hspace{0.01em}$

% (+++*** Orig. No.:  Definition Tree-TC )

\label{Definition Tree-TC}

If $t$ is a tree with root $<a,b>,$ then t/c will be the same tree,
only with the root $<c,b>.$

% karl-search= End Definition Tree-TC
\vspace{7mm}

% *************************************

\vspace{7mm}

% {\LARGE karl-search= Start Construction Smooth-Tree }

\index{Construction Smooth-Tree}

\ed

\bcs

$\hspace{0.01em}$

% (+++*** Orig. No.:  Construction Smooth-Tree )

\label{Construction Smooth-Tree}

(A) The set $T_{x}$ of trees $t$ for fixed $x$:

(1) Construction of the set $T \xbm_{x}$ of trees for those sets $U \xbe
\xdy,$ where $x \xbe \xbm (U):$

Let $U \xbe \xdy,$ $x \xbe \xbm (U).$ The trees $t_{U,x} \xbe T \xbm_{x}$
are constructed inductively,
observing simultaneously:

If $<U_{n+1},x_{n+1}>$ is a child of $<U_{n},x_{n}>,$ then
(a) $x_{n+1} \xbe \xbm (U_{n+1})-H(U_{n}),$
and
(b) $U_{n} \xcc U_{n+1}$.

Set $U_{0}:=U,$ $x_{0}:=x.$

Level 0: $<U_{0},x_{0}>.$

Level $n \xcp n+1$:
Let $<U_{n},x_{n}>$ be in level $n.$
Suppose $Y_{n+1} \xbe \xdy,$ $x_{n} \xbe Y_{n+1},$ and $Y_{n+1} \xcC
H(U_{n}).$ Note that $ \xbm (U_{n} \xcv Y_{n+1})-H(U_{n}) \xEd \xCQ $ by
Fact \ref{Fact HU-2}, (7), and $ \xbm (U_{n} \xcv Y_{n+1})-H(U_{n})
\xcc \xbm (Y_{n+1})$
by Fact \ref{Fact HU-2}, (3).
Choose $f_{n+1} \xbe \xbP \{ \xbm (U_{n} \xcv Y_{n+1})-H(U_{n}):$ $Y_{n+1}
\xbe \xdy,$ $x_{n} \xbe Y_{n+1} \xcC H(U_{n})\}$ (for the construction
of this tree, at this element), and let the set of
children of $<U_{n},x_{n}>$ be $\{<U_{n} \xcv Y_{n+1},f_{n+1}(Y_{n+1})>:$
$Y_{n+1} \xbe \xdy,$ $x_{n} \xbe Y_{n+1} \xcC H(U_{n})\}.$
(If there is no such $Y_{n+1}$, $<U_{n},x_{n}>$ has no children.)
Obviously, (a) and (b) hold.

We call such trees $U,x-$trees.

(2) Construction of the set $T'_{x}$ of trees for the nonminimal elements.
Let $x \xbe Z.$ Construct the tree $t_{x}$ as follows (here, one tree per
$x$ suffices for
all $U)$:

Level 0: $< \xCQ,x>$

Level 1:
Choose arbitrary $f \xbe \xbP \{ \xbm (U):x \xbe U \xbe \xdy \}.$ Note
that $U \xEd \xCQ \xcp \xbm (U) \xEd \xCQ $ by $Z=K$ (by
Remark \ref{Remark Gamma-x}, (1)). Let $\{<U,f(U)>:x \xbe U \xbe \xdy
\}$ be
the set of children of $< \xCQ,x>.$
This assures that the element will be nonminimal.

Level $>1$:
Let $<U,f(U)>$ be an element of level 1, as $f(U) \xbe \xbm (U),$ there is
a $t_{U,f(U)} \xbe T \xbm_{f(U)}.$
Graft one of these trees $t_{U,f(U)} \xbe T \xbm_{f(U)}$ at $<U,f(U)>$ on
the level 1.
This assures that a minimal element will be below it to guarantee
smoothness.

Finally, let $T_{x}:=T \xbm_{x} \xcv T'_{x}.$

(B) The relation $ \xej $ between trees:
For $x,y \xbe Z,$ $t \xbe T_{x}$, $t' \xbe T_{y}$, set $t \xem t' $ iff
for some $Y$ $<Y,y>$ is a child of
the root $<X,x>$ in $t,$ and $t' $ is the subtree of $t$ beginning at this
$<Y,y>.$

(C) The structure $ \xdz $:
Let $ \xdz $ $:=$ $<$ $\{<x,t_{x}>:$ $x \xbe Z,$ $t_{x} \xbe T_{x}\}$,
$<x,t_{x}> \xee <y,t_{y}>$ iff $t_{x} \xem^{*}t_{y}$ $>.$

% karl-search= End Construction Smooth-Tree
\vspace{7mm}

% *************************************

\vspace{7mm}

% {\LARGE karl-search= Start Fact Smooth-Tree }

\index{Fact Smooth-Tree}

\ecs

The rest of the proof are simple observations.

\bfa

$\hspace{0.01em}$

% (+++*** Orig. No.:  Fact Smooth-Tree )

\label{Fact Smooth-Tree}

(1) If $t_{U,x}$ is an $U,x-$tree, $<U_{n},x_{n}>$ an element of $t_{U,x}$
, $<U_{m},x_{m}>$ a direct or
indirect child of $<U_{n},x_{n}>,$ then $x_{m} \xce H(U_{n}).$

(2) Let $<Y_{n},y_{n}>$ be an element in $t_{U,x} \xbe T \xbm_{x}$, $t' $
the subtree starting at $<Y_{n},y_{n}>,$
then $t' $ is a $Y_{n},y_{n}-tree.$

(3) $ \xeb $ is free from cycles.

(4) If $t_{U,x}$ is an $U,x-$tree, then $<x,t_{U,x}>$ is $ \xeb -$minimal
in $ \xdz \xex U.$

(5) No $<x,t_{x}>,$ $t_{x} \xbe T'_{x}$ is minimal in any $ \xdz \xex U,$
$U \xbe \xdy.$

(6) Smoothness is respected for the elements of the form $<x,t_{U,x}>.$

(7) Smoothness is respected for the elements of the form $<x,t_{x}>$ with
$t_{x} \xbe T'_{x}.$

(8) $ \xbm = \xbm_{ \xdz }.$

% karl-search= End Fact Smooth-Tree
\vspace{7mm}

% *************************************

\vspace{7mm}

% {\LARGE karl-search= Start Fact Smooth-Tree Proof }

\index{Fact Smooth-Tree Proof}

\efa

\subparagraph{
Proof:
}

$\hspace{0.01em}$

% (+++*** Orig.:  Proof: )

\label{Section Proof:}

(1) trivial by (a) and (b).

(2) trivial by (a).

(3) Note that no $<x,t_{x}>$ $t_{x} \xbe T'_{x}$ can be smaller than any
other element (smaller
elements require $U \xEd \xCQ $ at the root). So no cycle involves any
such $<x,t_{x}>.$
Consider now $<x,t_{U,x}>,$ $t_{U,x} \xbe T \xbm_{x}$. For any
$<y,t_{V,y}> \xeb <x,t_{U,x}>,$ $y \xce H(U)$ by (1),
but $x \xbe \xbm (U) \xcc H(U),$ so $x \xEd y.$

(4) This is trivial by (1).

(5) Let $x \xbe U \xbe \xdy,$ then $f$ as used in the construction of
level 1 of $t_{x}$ chooses
$y \xbe \xbm (U) \xEd \xCQ,$ and some $<y,t_{U,y}>$ is in $ \xdz \xex U$
and below $<x,t_{x}>.$

(6) Let $x \xbe A \xbe \xdy,$ we have to show that either $<x,t_{U,x}>$
is minimal in $ \xdz \xex A,$ or that
there is $<y,t_{y}> \xeb <x,t_{U,x}>$ minimal in $ \xdz \xex A.$
Case 1, $A \xcc H(U)$: Then $<x,t_{U,x}>$ is minimal in $ \xdz \xex A,$
again by (1).
Case 2, $A \xcC H(U)$: Then A is one of the $Y_{1}$ considered for level
1. So there is
$<U \xcv A,f_{1}(A)>$ in level 1 with $f_{1}(A) \xbe \xbm (A) \xcc A$ by
Fact \ref{Fact HU-2}, (3).
But note that by (1)
all elements below $<U \xcv A,f_{1}(A)>$ avoid $H(U \xcv A).$ Let $t$ be
the subtree of $t_{U,x}$
beginning at $<U \xcv A,f_{1}(A)>,$ then by (2) $t$ is one of the $U \xcv
A,f_{1}(A)-trees,$ and
$<f_{1}(A),t>$ is minimal in $ \xdz \xex U \xcv A$ by (4), so in $ \xdz
\xex A,$ and $<f_{1}(A),t> \xeb <x,t_{U,x}>.$

(7) Let $x \xbe A \xbe \xdy,$ $<x,t_{x}>,$ $t_{x} \xbe T_{x}',$ and
consider the subtree $t$ beginning at $<A,f(A)>,$
then $t$ is one of the $A,f(A)-$trees, and $<f(A),t>$ is minimal in $ \xdz
\xex A$ by (4).

(8) Let $x \xbe \xbm (U).$ Then any $<x,t_{U,x}>$ is $ \xeb -$minimal in $
\xdz \xex U$ by (4), so $x \xbe \xbm_{ \xdz }(U).$
Conversely, let $x \xbe U- \xbm (U).$ By (5), no $<x,t_{x}>$ is minimal in
$U.$ Consider now some
$<x,t_{V,x}> \xbe \xdz,$ so $x \xbe \xbm (V).$ As $x \xbe U- \xbm (U),$
$U \xcC H(V)$ by Fact \ref{Fact HU-2}, (6).
Thus $U$ was
considered in the construction of level 1 of $t_{V,x}.$ Let $t$ be the
subtree of $t_{V,x}$
beginning at $<V \xcv U,f_{1}(U)>,$ by $ \xbm (V \xcv U)-H(V) \xcc \xbm
(U)$ (Fact \ref{Fact HU-2}, (3)),
$f_{1}(U) \xbe \xbm (U) \xcc U,$ and $<f_{1}(U),t> \xeb <x,t_{V,x}>.$

$ \xcz $
\\[3ex]

% karl-search= End Fact Smooth-Tree Proof
\vspace{7mm}

% *************************************

\vspace{7mm}

% {\LARGE karl-search= Start Proposition A-Smooth-Trans-Complete Proof }

\index{Proposition A-Smooth-Trans-Complete Proof}

\paragraph{
Proof of Proposition \ref{Proposition A-Smooth-Trans-Complete}
}

$\hspace{0.01em}$

% (+++*** Orig.:  Proof of  {Proposition A-Smooth-Trans-Complete } )

\label{Section Proof of  {Proposition A-Smooth-Trans-Complete }}

Consider the construction in the proof of
Proposition \ref{Proposition Smooth-Complete-Trans}.

Thus, we only have to show that in $ \xdz $ defined by

$ \xdz $ $:=$ $<$ $\{<x,t_{x}>:$ $x \xbe Z,$ $t_{x} \xbe T_{x}\}$,
$<x,t_{x}> \xee <y,t_{y}>$ iff $t_{x} \xem^{*}t_{y}$ $>,$ $t_{x} \xem
t_{y}$ implies
$rg(y) \xck rg(x).$

But by construction of the trees,
$x_{n} \xbe Y_{n+1},$ and $x_{n+1} \xbe \xbm (U_{n} \xcv Y_{n+1}),$ so
$rg(x_{n+1}) \xck rg(x_{n}).$

$ \xcz $ (Proposition \ref{Proposition A-Smooth-Trans-Complete})
\\[3ex]

% karl-search= End Proposition A-Smooth-Trans-Complete Proof
\vspace{7mm}

% *************************************

\vspace{7mm}

% karl-search= End Hier-ARepr-Smooth-Trans
\vspace{7mm}

% *************************************

\vspace{7mm}

% karl-search= End Hier-ARepr-Smooth
\vspace{7mm}

% *************************************

\vspace{7mm}

%  6.5  The logical properties with dp
%  6.5  The logical properties with dp
% %
% ====================================
\subsection{
The logical properties with definability preservation
}

% {\LARGE karl-search= Start Hier-ARepr-Logic }

\label{Section Hier-ARepr-Logic}
\index{Section Hier-ARepr-Logic}

First, a small fact about the $ \xda.$

\bfa

$\hspace{0.01em}$

% (+++*** Orig. No.:  Fact 6.21 )

\label{Fact 6.21}

Let $ \xda $ be as above (and thus finite).
Then each $A_{i}$ is equivalent to a formula $ \xba_{i}.$

\efa

\subparagraph{
Proof:
}

$\hspace{0.01em}$

% (+++*** Orig.:  Proof: )

\label{Section Proof:}

We use the standard topology and its compactness.
By definition, each $M(A_{i})$ is closed, by finiteness all unions of such
$M(A_{i})$ are
closed, too, so $ \xdC (M(A_{i}))$ is closed. By compactness, each open
cover $X_{j}:j \xbe J$
of the clopen $M(A_{i})$ contains a finite subcover, so also $ \xcV
\{M(A_{j}):j \xEd i\}$ has
a finite open cover. But the $M( \xbf ),$ $ \xbf $ a formula form a basis
of the closed
sets, so we are done. $ \xcz $
\\[3ex]

\bp

$\hspace{0.01em}$

% (+++*** Orig. No.:  Proposition 6.22 )

\label{Proposition 6.22}

Let $ \xcn $ be a logic for $ \xdl.$ Set $T^{ \xdm }:=Th( \xbm_{ \xdm
}(M(T))),$
where $ \xdm $ is a preferential structure.

(1) Then there is a (transitive) definability preserving
classical preferential model $ \xdm $ s.t. $ \ol{ \ol{T} }=T^{ \xdm }$ iff

(LLE), (CCL), (SC), (PR) hold for all $T,T' \xcc \xdl.$

(2) The structure can be chosen smooth, iff, in addition

(CUM) holds.

(3) The structure can be chosen $ \xda -$ranked, iff, in addition

$( \xda -$min) $T \xcL \xCN \xba_{i}$ and $T \xcL \xCN \xba_{j},$ $i<j$
implies $ \ol{ \ol{T} } \xcl \xCN \xba_{j}$

holds.

\ep

The proof is an immediate consequence of Proposition \ref{Proposition 6.23}
and the respective above results.
This proposition (or its analogue) was mostly already
shown in  \cite{Sch92} and
 \cite{Sch96-1} and is repeated here for completeness' sake, but
with a new and partly stronger proof.

\bp

$\hspace{0.01em}$

% (+++*** Orig. No.:  Proposition 6.23 )

\label{Proposition 6.23}

Consider for a logic $ \xcn $ on $ \xdl $ the properties

(LLE), (CCL), (SC), (PR), (CUM), $( \xda -$min) hold for all $T,T' \xcc
\xdl.$

and for a function $ \xbm: \xdD_{ \xdl } \xcp \xdp (M_{ \xdl })$ the
properties

$( \xbm dp)$ $ \xbm $ is definability preserving, i.e. $ \xbm (M(T))=M(T'
)$ for some $T' $

$( \xbm \xcc ),$ $( \xbm PR),$ $( \xbm CUM),$ $( \xbm \xda )$

for all $X,Y \xbe \xdD_{ \xdl }.$

It then holds:

(a) If $ \xbm $ satisfies $( \xbm dp),$ $( \xbm \xcc ),$ $( \xbm PR),$
then $ \xcn $ defined by $ \ol{ \ol{T} }:=T^{ \xbm }:=$
$Th( \xbm (M(T)))$ satisfies (LLE), (CCL), (SC), (PR).
If $ \xbm $ satisfies in addition $( \xbm CUM),$ then (CUM) will hold,
too.
If $ \xbm $ satisfies in addition $( \xbm \xda ),$ then $( \xda -$min)
will hold, too.

(b) If $ \xcn $ satisfies (LLE), (CCL), (SC), (PR),
then there is $ \xbm: \xdD_{ \xdl } \xcp \xdp (M_{ \xdl })$ s.t. $ \ol{
\ol{T} }=T^{ \xbm }$
for all $T \xcc \xdl $ and $ \xbm $ satisfies $( \xbm dp),$ $( \xbm \xcc
),$ $( \xbm PR).$
If, in addition, (CUM) holds, then $( \xbm CUM)$ will hold, too.
If, in addition, $( \xda -$min) holds, then $( \xbm \xda )$ will hold,
too.

\ep

\subparagraph{
Proof:
}

$\hspace{0.01em}$

% (+++*** Orig.:  Proof: )

\label{Section Proof:}

Set $ \xbm (T):= \xbm (M(T)),$ note that $ \xbm (T \xcv T' ):= \xbm (M(T
\xcv T' ))= \xbm (M(T) \xcs M(T' )).$

(a)
Suppose $ \ol{ \ol{T} }=T^{ \xbm }for$ some such $ \xbm,$ and all $T.$

(LLE): If $ \ol{T}= \ol{T' }$, then $M(T)=M(T' ),$ so $ \xbm (T)= \xbm
(T' ),$ and $T^{ \xbm }=T'^{ \xbm }$.

(CCL) and (SC) are trivial.

We show (PR):
$M( \ol{ \ol{T} } \xcv T' )=M( \ol{ \ol{T} }) \xcs M(T' )=_{( \xbm dp)}
\xbm (T) \xcs M(T' )=_{( \xbm \xcc )} \xbm (T) \xcs M(T) \xcs M(T' )=$
$ \xbm (T) \xcs M(T \xcv T' ) \xcc_{( \xbm PR)} \xbm (T \xcv T' )=_{( \xbm
dp)}M( \ol{ \ol{T \xcv T' } }).$
Let now $ \xbf \xbe \ol{ \ol{T \xcv T' } }$, so $ \xbf $ holds in all $m
\xbe M( \ol{ \ol{T \xcv T' } }),$ so $ \xbf $ holds in
all $m \xbe M( \ol{ \ol{T} } \xcv T' ),$ so $ \ol{ \ol{T} } \xcv T' \xcl
\xbf,$ so $ \xbf \xbe \ol{ \ol{ \ol{T} } \xcv T' }.$

We turn to (CUM):

Let $T$ $ \xcc $ $ \ol{T' }$ $ \xcc $ $ \ol{ \ol{T} }.$ Thus by $( \xbm
Cum)$ and $ \xbm (T)$ $ \xcc $ $M( \ol{ \ol{T} })$ $ \xcc $ $M(T' )$ $
\xcc $ $M(T),$ so $ \xbm (T)= \xbm (T' ),$
so $ \ol{ \ol{T} }$ $=$ $Th( \xbm (T))$ $=$ $Th( \xbm (T' ))$ $=$ $ \ol{
\ol{T' } }.$

$( \xda -$min) is trivial.

(b)
Let $ \xcn $ satisfy $(LLE)-(CUM)$ for all $T.$ We define $ \xbm $ and
show $ \ol{ \ol{T} }=T^{ \xbm }.$ (CUM) will be
needed only to show $( \xbm CUM).$

If $X=M(T)$ for some $T \xcc \xdl,$ set $ \xbm (X):=M( \ol{ \ol{T} }).$

If $X=M(T)=M(T' )$, then $ \ol{T}= \ol{T' }$, thus $ \ol{ \ol{T} }= \ol{
\ol{T' } }$ by (LLE), so $M( \ol{ \ol{T} })$ $=M( \ol{ \ol{T' } }),$
and $ \xbm $ is well-defined.
Moreover, $ \xbm $ satisfies $( \xbm dp),$ and by (SC), $ \xbm (X) \xcc X$
. We show $ \ol{ \ol{T} }=T^{ \xbm }$:
Let now $T \xcc \xdl $ be given. Then $ \xbf \xbe T^{ \xbm }$ $: \xcr $ $
\xcA m \xbe \xbm (T).m \xcm \xbf $ $ \xcr $
$ \xcA m \xbe M( \ol{ \ol{T} }).m \xcm \xbf $ $ \xcr $ $ \ol{ \ol{T} }
\xcl \xbf $ $ \xcr $ $ \xbf \xbe \ol{ \ol{T} }$ (as $ \ol{ \ol{T} }$ is
classically closed).

Next, we show that the above defined $ \xbm $ satisfies $( \xbm PR).$
Suppose $X:=M(T),$ $Y:=M(T' ),$ $X \xcc Y,$ we have to show $ \xbm (Y)
\xcs X \xcc \xbm (X).$
By prerequisite, $ \ol{T' } \xcc \ol{T},$ so $ \ol{T \xcv T' }= \ol{T},$
so $ \ol{ \ol{T \xcv T' } }= \ol{ \ol{T} }$ by (LLE). By (PR) $ \ol{ \ol{T
\xcv T' } } \xcc \ol{ \ol{ \ol{T' } } \xcv T},$
so $ \xbm (Y) \xcs X= \xbm (T' ) \xcs M(T)=M( \ol{ \ol{T' } } \xcv T) \xcc
M( \ol{ \ol{T \xcv T' } })=M( \ol{ \ol{T} })= \xbm (X),$ using $( \xbm
dp).$

$( \xbm \xda )$ is trivial.

It remains to show $( \xbm CUM).$
So let $X=M(T),$ $Y=M(T' ),$ and $ \xbm (T):=M( \ol{ \ol{T} }) \xcc M(T' )
\xcc M(T)$ $ \xcp $ $ \ol{T} \xcc \ol{T' } \xcc \ol{ \ol{T} }= \ol{ \ol{(
\ol{T})} }$ $ \xcp $ $ \ol{ \ol{T} }= \ol{ \ol{(
\ol{T})} }= \ol{ \ol{( \ol{T' })} }= \ol{ \ol{T' } }$ $ \xcp $
$ \xbm (T)=M( \ol{ \ol{T} })=M( \ol{ \ol{T' } })= \xbm (T' ),$ thus $ \xbm
(X)= \xbm (Y).$

$ \xcz $ (Proposition \ref{Proposition 6.23})
\\[3ex]

% karl-search= End Hier-ARepr-Logic
\vspace{7mm}

% *************************************

\vspace{7mm}

% karl-search= End Hier-ARepr
\vspace{7mm}

% *************************************

\vspace{7mm}

%  7  Formal results and representation for hierarchical conditionals
%  7  Formal results and representation for hierarchical conditionals
% %
% =============================
\section{
Formal results and representation for hierarchical conditionals
}

% {\LARGE karl-search= Start Hier-CondRepr }

\label{Section Hier-CondRepr}
\index{Section Hier-CondRepr}

We look here at the following problem:

Given

(1.1) a finite, ordered partition $ \xda $ of $ \xdA,$ $ \xda =<\{A_{i}:i
\xbe I\},<>$

(1.2) a normality relation $ \xeb,$ which is an $ \xda -$ranking,
defining a
choice function $ \xbm $ on subsets of $ \xdA,$ (so, obviously,
$A<A' $ iff $ \xbm (A \xcv A' ) \xcs A' = \xCQ ),$

(1.3) a subset $ \xdB \xcc \xdA,$ and we set $ \xdc:=< \xda, \xdB >$
(thus, the $B_{i}$ are just $A_{i} \xcs \xdB,$
this way of writing saves a little notation),

(2.1) a set of models $M,$

(2.2) an accessibility relation $R$ on $M,$ with some finite upper bound
on
$R-$chains,

(2.3) an unknown extension of $R$ to pairs $(m,a),$ $m \xbe M,$ $a \xbe
\xdA,$

(3.1) a notion of validity $m \xcm \xdc,$ for $m \xbe M,$ defined by
$m \xcm \xdc $ iff
$\{a \xbe \xdA:$ $ \xcE A \xbe \xda (a \xbe \xbm (A),$ $a \xbe R(m),$ and

$ \xDC \xDC \xDC \xCN \xcE a' ( \xcE A' \xbe \xda (a' \xbe \xbm (A' ),$
$a' \xbe R(m),$ $a' \xeb a\}$ $ \xcc $ $ \xdB,$

(3.2) a subset $M' $ of $M$

give a criterion which decides whether it is possible to construct the
extension of $R$ to pairs $(m,a)$ s.t. $ \xcA m \xbe M.(m \xbe M' $ $ \xcj
$ $m \xcm \xdc ).$

We first show some elementary facts on the situation, and give the
criterion
in Proposition \ref{Proposition 7.4}, together with the proof that it does
what is wanted.

\bfa

$\hspace{0.01em}$

% (+++*** Orig. No.:  Fact 7.1 )

\label{Fact 7.1}

Reachability for a transitive relation is characterized by

$y \xbe R(x)$ $ \xcp $ $R(y) \xcc R(x)$

\efa

\subparagraph{
Proof:
}

$\hspace{0.01em}$

% (+++*** Orig.:  Proof: )

\label{Section Proof:}

Define directly xRz iff $z \xbe R(x).$ This does it. $ \xcz $
\\[3ex]

Let now $S$ be a set with an accessibility relation $R',$ generated by
transitive
closure from the intransitive subrelation $R.$ All modal notation will be
relative
to this $R.$

Let $ \xdA =M( \xba ),$ $A_{i}=M( \xba_{i}),$ the latter is justified by
Fact \ref{Fact 6.21}.

\bd

$\hspace{0.01em}$

% (+++*** Orig. No.:  Definition 7.1 )

\label{Definition 7.1}

(1) $ \xbm ( \xda ):= \xcV \{ \xbm (A_{i}):i \xbe I\}$

(warning: this is NOT $ \xbm ( \xdA ))$

(2) $ \xda_{m}:=R(m) \xcs \xdA,$

(3) $ \xbm ( \xda_{m}):= \xcV \{ \xbm (A_{i}) \xcs R(m):i \xbe I\}$

(3a) $ \xbn ( \xda_{m})$ $:=$ $ \xbm ( \xbm ( \xda_{m}))$

(thus $ \xbn ( \xda_{m})$ $=$ $\{a \xbe \xdA:$ $ \xcE A \xbe \xda (a \xbe
\xbm (A),$ $a \xbe R(m),$ and

$ \xDC \xDC \xDC \xCN \xcE a' ( \xcE A' \xbe \xda (a' \xbe \xbm (A' ),$
$a' \xbe R(m),$ $a' \xeb a\}.$

(4) $m \xcm \xdc $ $: \xcr $ $ \xbn ( \xda_{m})) \xcc \xdB.$

See Diagram \ref{Diagram C-Validity}

\vspace{20mm}

\begin{diagram}

\label{Diagram C-Validity}
\index{Diagram C-Validity}

\centering
\setlength{\unitlength}{0.00083333in}
{\renewcommand{\dashlinestretch}{30}
\begin{picture}(3841,3665)(0,200)
\put(2407.000,2628.000){\arc{1110.000}{3.4719}{5.9529}}
\put(2407.000,1578.000){\arc{1110.000}{3.4719}{5.9529}}
\put(2398.149,661.355){\arc{1503.312}{3.9463}{5.4292}}
\put(-24.763,2227.012){\arc{4826.591}{5.7407}{6.7841}}
\put(2369,2283){\ellipse{1874}{2550}}
\path(1507,2808)(3232,2808)
\path(1507,1758)(3232,1758)
\path(907,3263)(2562,1963)
\path(2449.100,2013.534)(2562.000,1963.000)(2486.163,2060.718)
\path(157,3248)(2072,1213)
\path(1967.915,1279.831)(2072.000,1213.000)(2011.611,1320.950)
\path(272,3258)(847,3258)
\put(210,3258){\circle*{20}}
\put(880,3258){\circle*{20}}
\path(727.000,3228.000)(847.000,3258.000)(727.000,3288.000)
\put(3457,3108){{\xssc $A''$}}
\put(3457,2208){{\xssc $A'$}}
\put(3457,1308){{\xssc $A$}}
\put(4000,1093){{\xssc $m \xcm \xdc$}}
\put(4000,8033){{\xssc $m' \xcM \xdc$}}
\put(2407,2883){{\xssc $\mu(A'')$}}
\put(2407,1833){{\xssc $\mu(A')$}}
\put(2407,1093){{\xssc $\mu(A)$}}
\put(132,3318){{\xssc $m$}}
\put(882,3313){{\xssc $m'$}}
\put(1972,3518){{\xssc $B$}}
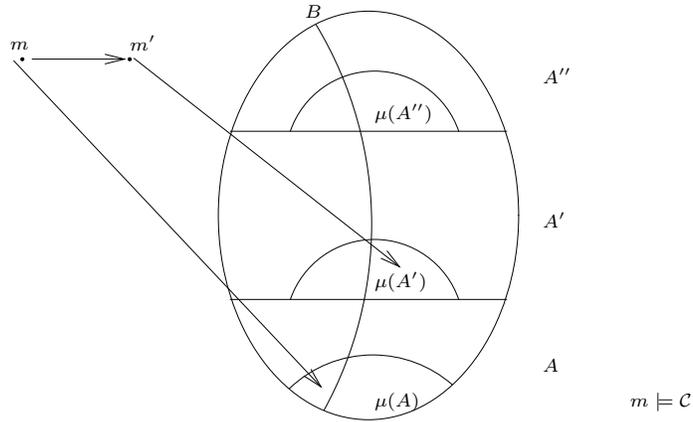
\put(10,800){{\xssc Here, the ``best'' element $m$ sees is in
$B$, so $\xdc$ holds in $m$.}}
\put(10,600){{\xssc The ``best'' element $m'$ sees is not in
$B$, so $\xdc$ does not hold in $m'$.}}

\put(100,250) {{\rm\bf Validity of $\xdc$ from $m$ and $m'$}}

\end{picture}
}
\end{diagram}

\vspace{4mm}

\ed

We have the following Fact for $m \xcm \xdc:$

\bfa

$\hspace{0.01em}$

% (+++*** Orig. No.:  Fact 7.2 )

\label{Fact 7.2}

Let $m,m' \xbe M,$ $A,A' \xbe \xda.$

(1) $m \xcm \xcX \xCN \xba,$ mRm' $ \xch $ $m' \xcm \xcX \xCN \xba $

(2) $ \xCf mRm',$ $ \xbn ( \xda_{m}) \xcs A \xEd \xCQ,$ $ \xbn (
\xda_{m' }) \xcs A' \xEd \xCQ,$ $ \xch $ $A \xck A' $

(3) $ \xCf mRm',$ $ \xbn ( \xda_{m}) \xcs A \xEd \xCQ,$ $ \xbn (
\xda_{m' }) \xcs A' \xEd \xCQ,$ $m \xcm \xdc,$ $m' \xcM \xdc, \xch $
$A<A' $

\efa

\subparagraph{
Proof:
}

$\hspace{0.01em}$

% (+++*** Orig.:  Proof: )

\label{Section Proof:}

Trivial. $ \xcz $
\\[3ex]

\bfa

$\hspace{0.01em}$

% (+++*** Orig. No.:  Fact 7.3 )

\label{Fact 7.3}

We can conclude from above properties that there are no arbitrarily long
$R-$chains of models $m,$ changing from $m \xcm \xdc $ to $m \xcM \xdc $
and back.

\efa

\subparagraph{
Proof:
}

$\hspace{0.01em}$

% (+++*** Orig.:  Proof: )

\label{Section Proof:}

Trivial: By Fact \ref{Fact 7.2}, (3), any change from $ \xcm \xdc $
to $ \xcM \xdc $ results in a
strict increase in rank. $ \xcz $
\\[3ex]

We solve now the representation task described at the beginning of
Section \ref{Section Hier-CondRepr}, all we need
are the properties shown in Fact \ref{Fact 7.2}.

(Note that constructing $R$ between the different $m,$ $m' $ is trivial:
we could just
choose the empty relation.)

\bp

$\hspace{0.01em}$

% (+++*** Orig. No.:  Proposition 7.4 )

\label{Proposition 7.4}

If the properties of Fact \ref{Fact 7.2} hold, we can extend $R$ to
solve the
representation problem described at the beginning of
this Section \ref{Section Hier-CondRepr}.

\ep

\subparagraph{
Proof:
}

$\hspace{0.01em}$

% (+++*** Orig.:  Proof: )

\label{Section Proof:}

By induction on $R.$ This is possible, as the depth of $R$ on $M$ was
assumed to
be finite.

\bcs

$\hspace{0.01em}$

% (+++*** Orig. No.:  Construction 7.1 )

\label{Construction 7.1}

We choose now elements as possible, which ones are chosen exactly does not
matter.

$X_{i}:=\{b_{i},c_{i}\}$ iff $ \xbm (A_{i}) \xcs \xdB \xEd \xCQ $ and $
\xbm (A_{i})- \xdB \xEd \xCQ,$ $b_{i} \xbe \xbm (A_{i}) \xcs \xdB,$
$c_{i} \xbe \xbm (A_{i})- \xdB.$

$X_{i}:=\{c_{i}\}$ iff $ \xbm (A_{i}) \xcs \xdB = \xCQ $ and $ \xbm
(A_{i})- \xdB \xEd \xCQ,$ $c_{i} \xbe \xbm (A_{i})- \xdB $

$X_{i}:=\{b_{i}\}$ iff $ \xbm (A_{i}) \xcs \xdB \xEd \xCQ $ and $ \xbm
(A_{i})- \xdB = \xCQ,$ $b_{i} \xbe \xbm (A_{i}) \xcs \xdB,$

$X_{i}:= \xCQ $ iff $ \xbm (A_{i})= \xCQ.$

Case 1:

Let $m$ be $R-$minimal and $m \xcm \xdc.$ Let $i_{0}$ be the first $i$
s.t. $b_{i} \xbe X_{i},$ make
$ \xbg (m):=i_{0},$ and make $R(m):=\{b_{i_{0}}\} \xcv \xcV
\{X_{i}:i>i_{0}\}.$ This makes $ \xdc $ hold.
(This leaves us as many possibilities open as possible - remember we have
to
decrease the set of reachable elements now.)

Case 2:

Let $m$ be $R-$minimal and $m \xcM \xdc.$ Let $i_{0}$ be the first $i$
s.t. $c_{i} \xbe X_{i},$ make
$ \xbg (m):=i_{0},$ and make $R(m):= \xcV \{X_{i}:i \xcg i_{0}\}.$ This
makes $ \xdc $ false.

Let all $R-$predecessors of $m$ be determined, and $i:=max\{ \xbg (m' ):m'
Rm\}.$

Case 3: $m \xcm \xdc.$ Let $j$ be the smallest $i' \xcg i$ with $ \xbm
(A_{i' }) \xcs \xdB \xEd \xCQ.$
Let $R(m):=\{b_{j}\} \xcv \xcV \{X_{k}:k>j\},$ and $ \xbg (m):=j.$

Case 4: $m \xcM \xdc.$

Case 4.1: For all $m' Rm$ with $i= \xbg (m' )$ $m' \xcM \xdc.$

Take one such $m' $ and set $R(m):=R(m' ),$ $ \xbg (m):=i.$

Case 4.2: There is $m' Rm$ with $i= \xbb (m' )$ $m' \xcm \xdc.$

Let $j$ be the smallest $i' >i$ with $ \xbm (A_{i' })- \xdB \xEd \xCQ.$
Let $R(m):= \xcV \{X_{k}:k \xcg j\}.$
(Remark: To go from $ \xcm $ to $ \xcM,$ we have to go higher in the
hierarchy.)

Obviously, validity is done as it should be.
It remains to show that the sets of reachable elements decrease with $R.$

\ecs

\bfa

$\hspace{0.01em}$

% (+++*** Orig. No.:  Fact 7.5 )

\label{Fact 7.5}

In above construction, if mRm', then $R(m' ) \xcc R(m).$

\efa

\subparagraph{
Proof:
}

$\hspace{0.01em}$

% (+++*** Orig.:  Proof: )

\label{Section Proof:}

By induction, considering $R.$ $ \xcz $ (Fact \ref{Fact 7.5} and
Proposition \ref{Proposition 7.4})
\\[3ex]

We consider an example for illustration.

\be

$\hspace{0.01em}$

% (+++*** Orig. No.:  Example 7.1 )

\label{Example 7.1}

Let $a_{1}Ra_{2}RcRc_{1},$ $b_{1}Rb_{2}Rb_{3}RcRd_{1}Rd_{2}.$

Let $ \xdc =(A_{1}>B_{1}, \Xl,A_{n}>B_{n})$ with the $C_{i}$ consistency
with $ \xbm (A_{i}).$

Let $ \xbm (A_{2}) \xcs B_{2}= \xCQ,$ $ \xbm (A_{3}) \xcc B_{3},$ and for
the other $i$ hold neither of these two.

Let $a_{1},a_{2},b_{2},c_{1},d_{2} \xcm \xdc,$ the others $ \xcM \xdc.$

Let $ \xbm (A_{1})=\{a_{1,1},a_{1,2}\},$ with $a_{1,1} \xbe B_{1},$
$a_{1,2} \xce B_{1},$

$ \xbm (A_{2})=\{a_{2,1}\},$ $ \xbm (A_{3})=\{a_{3,1}\}$ (there is no
reason to differentiate),

and the others like $ \xbm (A_{1}).$ Let $ \xbm A:= \xcV \{ \xbm (A_{i}):i
\xck n\}.$

We have to start at $a_{1}$ and $b_{1},$ and make $R(x)$ progressively
smaller.

Let $R(a_{1}):= \xbm A-\{a_{1,2}\},$ so $a_{1} \xcm \xdc.$ Let
$R(a_{2})=R(a_{1}),$ so again $a_{2} \xcm \xdc.$

Let $R(b_{1}):= \xbm A-\{a_{1,1}\},$ so $b_{1} \xcM \xdc.$ We now have to
take $a_{1,2}$ away, but
$a_{2,1}$ too to be able to change. So let
$R(b_{2}):=R(b_{1})-\{a_{1,2},a_{2,1}\},$ so we
begin at $ \xbm (A_{3}),$ which is a (positive) singleton.
Then let $R(b_{3}):=R(b_{2})-\{a_{3,1}\}.$

We can choose $R(c):=R(b_{3}),$ as $R(b_{3}) \xcc R(a_{2}).$

Let $R(c_{1}):=R(c)-\{a_{4,2}\}$ to make $ \xdc $ hold again. Let
$R(d_{1}):=R(c),$ and $R(d_{2}):=R(c_{1}).$

$ \xcz $
\\[3ex]

% karl-search= End Hier-CondRepr
\vspace{7mm}

% *************************************

\vspace{7mm}

\ee

\end{document}